\documentclass[11pt]{article}
\usepackage{geometry}                % See geomet\rho Iy.pdf to lea\rho In the layout options. The\rho Ie a\rho Ie lots.
\usepackage{authblk}
\usepackage{graphicx}
\usepackage{amssymb}
\usepackage{amsthm}
\usepackage{amsmath}
\usepackage{mathabx}
\usepackage{enumitem}
\usepackage{silence}
\usepackage[square,numbers]{natbib}
\usepackage{xcolor}
\usepackage{algorithm}
\usepackage[noend]{algpseudocode}
\usepackage[section]{placeins}

\usepackage{titlesec}

\titleformat{\subsection}{\normalfont\normalsize}{\thesubsection}{1em}{}
\usepackage{hyperref}
\usepackage{nameref}
\usepackage{cleveref} 
\usepackage{url}

\usepackage{ulem}
\newlength{\maxwidth}
\newcommand{\algalign}[2]% #1 = text to left, #2 = text to right
{\makebox[\maxwidth][r]{$#1{}$}${}#2$}

\newtheorem{remark}{Remark}[section]
\crefname{remark}{Remark}{Remarks}

\newtheorem{assumption}{Assumption}[section]
\crefname{assumption}{Assumption}{Assumptions}

\newtheorem{definition}{Definition}[section]
\crefname{definition}{Definition}{Definitions}

\newtheorem{theorem}{Theorem}[section]
\crefname{theorem}{Theorem}{Theorems}

\newtheorem{claim}{Claim}[section]
\crefname{claim}{Claim}{claims}

\newtheorem{proposition}{{Proposition}}[section]
\crefname{proposition}{Proposition}{Propositions}

\newtheorem{lemma}{{Lemma}}[section]
\crefname{lemma}{Lemma}{Lemmas}

\newtheorem{corollary}{{Corollary}}[section]
\crefname{corollary}{Corollary}{Corollaries}

\crefname{algorithm}{algorithm}{algorithms}

\newtheorem{example}{{Example}}[section]
\crefname{example}{Example}{Examples}

\title{Hyperparameter Selection via Early Stopping for Bayesian Semilinear PDEs}

\author[1]{Maia Tienstra \footnote{Corresponding author: mtienstra@insa-toulouse.fr}}
\author[2]{Gottfried Hastermann}
\affil[1]{Institut de Mathématiques de Toulouse, Université de Toulouse, INSA}
\affil[2]{Institut für Mathematik, Technische Universität Ilmenau}

\begin{document}

%%%%%%%%%%%%%%%%%%%%%%%%%%%%%%%%%%%%%%
\maketitle
%%%%%%%%%%%%%%%%%%%%%%%%%%%%%%%%%%%%%%

%%%%%%%%%%%%%%%%%%%%%%%%%%%%%%%%%%%%%%
\begin{abstract}
We provide a data-driven method for choosing the scale of a Gaussian prior for non-linear Bayesian inverse problems arising from semilinear PDEs. Following \cite{koers2024}, the non-linear model for the parameter $f$ is reparametrized as a linear model for $v = \mathbb{L}u_f$, in which the prior scale can be selected by early stopping at the discrepancy principle \cite{tienstra2025}. We extend the transfer of frequentist guarantees from the linearised problem to the original one beyond Lipschitz solution maps, to maps admitting an arbitrary modulus of continuity, and we provide a general condition on the non-linearity under which such a map exists and is Lipschitz near the truth. The resulting posterior for $f$ contracts adaptively over a range of Sobolev smoothness, and its credible sets have asymptotic frequentist coverage one. We demonstrate our theory in detail for the stationary Schr\"odinger equation, for which we also provide numerical experiments. We further show how the results apply to the stationary Allen-Cahn equation and to a parabolic Allen-Cahn equation, and we discuss Darcy flow as an edge case outside the semilinear setting. The proposed method thus provides a data-driven way to tune Gaussian priors via early stopping, which is computationally efficient and statistically near-optimal for non-linear problems.
\end{abstract}

%%%%%%%%%%%%%%%%%%%%%%%%%%%%%%%%%%%%%%

\section{Introduction} \label{sec:intro}
 
%%%%%%%%%%%%%%%%%%%%%%%%%%%%%%%%%%%%%%
 
 Bayesian methods for parameter estimation of partial differential equations (PDEs) have emerged as an important field of research in the past decade, see \cite{Ghosal2017book}, \cite{koers2024}, \cite{Nickl2020BVM_Schro}, \cite{NicklLectureNotes2023}, \cite{NicklWang}, \cite{Stuart}, and references therein. A major motivation for applying Bayesian methods is the uncertainty quantification of the resulting point estimator. The literature has shown that prior choice plays a critical role in the performance of Bayesian methods in the non-parametric setting \cite{schwartz1965bayes} and \cite{Ghosal2017book}. To ensure optimal posterior contraction, consistency, and frequentist coverage, one must carefully choose the prior even in the linear (non-parametric) setting \cite{knapik2011}. In this paper, we are interested in studying early stopping as a prior selection method for non-linear Bayesian inverse problems which arise from PDEs. We focus on the problem of inferring a parameter $f$ of the underlying PDE with known boundary conditions and source function from observations
 \begin{align}
   \label{eq:obs}
    Y_i &:= G(f)(X_i) + \epsilon_i \\
     & = u_f(X_i) + \epsilon_i
\end{align}
 where $u_f$ is the solution to the $2m$ order semilinear partial differential equation 
\begin{equation}
\label{eq:generic_pde_full_operator}
     \left\{
    \begin{alignedat}{4}
        \mathcal{L}_{f}(u) &=& \, h & \text{ on } \mathcal{O}, \\
        \frac{\partial^{(i)}u}{\partial^{(i)} \nu}  &=& \, g_i  & \text{ on } \partial\mathcal{O} .
    \end{alignedat}
    \right.
\end{equation}
where $h$ and $g_i$ are known for every $i\in \{0,\dots, 2m-1\}$ and $ \epsilon_i \sim \mathcal{N}(0,1)$ and $\mathcal{O} \subset \mathbb{R}^d$. 
We suppose that we have $n$ such observations and will denote the collection of these observations by $Y_n$.

\begin{example}[Stationary Schrödinger Equation]
\label{example:schroedinger_eq}
The stationary Schrödinger equation is the guiding example for this work. The full treatment of this problem with the contraction rate and coverage, is given in \cref{corollary:schrodinger_contraction_rate}.
Let $\mathcal{O}$ be a bounded domain in $\mathbb{R}^d$ and let $f\in \mathcal{F}\subseteq L^\infty(\mathcal{O})$
Then the equation is given by
 
\begin{equation}
    \label{eq:schroedinger_eq}
    \left\{
    \begin{alignedat}{4}
        - \Delta u_f + fu &=& \, h & \text{ on } \mathcal{O}, \\
        u_f &=& \, g  & \text{ on } \partial \mathcal{O}.
    \end{alignedat}
    \right.
\end{equation}
\end{example}
 
To infer  $f$, we will use the Bayesian approach, which requires one to select a prior distribution for
$f$. Given this prior distribution and a forward operator $G$, we can derive a posterior distribution $\Pi(f \mid Y_n)$. The Bayesian method hence provides an entire distribution for $f$ conditional on observations \cref{eq:obs}. We can, in theory, compute a point estimator for $f$ from $\Pi(f \mid Y_n)$ by computing the mode. We will consider $G$ to be fixed and known, and thus what we can choose is the prior. The goal of this paper, then, is to select the best prior given a family of prior distributions for $f$, indexed by $\tau$ and denoted as $\Pi_\tau(f)$.

%%%%%%%%%%%%%%%%%%%%%%%%%%%%%%%%%%%%%%
\subsection{Statistical Model and Bayesian Approach}\label{subsec:bayesian_setup}
%%%%%%%%%%%%%%%%%%%%%%%%%%%%%%%%%%%%%%
The Bayesian paradigm to infer $f$ is to place prior $\Pi_n(f)$ over $f$ and assume that there exists a unique $f_0$ that generates data and is a solution to \cref{eq:generic_pde_full_operator}.  We can formulate estimating $f$ given $u_f$ as a regression problem in the following way, see \cite{NicklLectureNotes2023}. Suppose we want to estimate a function $f: \mathcal{O} \rightarrow \mathbb{R}$, where $\mathcal{O} \subset \mathbb{R}^d$ is a bounded open subset, from noisy observations of $u_f$ which is the solution of the partial differential equation \cref{eq:generic_pde_full_operator}. We then fix a parameter space $\mathcal{F} \subset L^2(\mathcal{O})$ which is measurable with respect to $\nu$, and whose elements are strictly positive, so that $G$ is well defined, and define forward maps
\begin{equation}
    f \mapsto G(f), \quad G:\mathcal{F} \rightarrow L^2(\mathcal{O})
\end{equation}
where $G$ is the solution map $(h,f) \mapsto u_f$ of \cref{eq:generic_pde_full_operator}. We drop the dependence on $h$, as we assume it is a fixed, known quantity. We assume we can take measurements of $G(f)$, which in practical applications consists of discrete measurements of $u_f$ over a finite set $X_1, ...,X_n$ of $\mathcal{O}$ plus noise, which we model as in \cref{eq:obs}, and where we further assume that $X_i \overset{\rm i.i.d.}{\sim} Uniform(\mathcal{O})$. We then collect our data as $\mathcal{D}^{(n)} = \{Y_i, X_i\}_{i=1}^n$, and the product measure of the joint law of the random variables $\mathcal{D}^{(n)}$ will be denoted as $P_f^n := \otimes_{i=1}^n P_f^i$. From the observations \cref{eq:obs}, the log-likelihood is 
\begin{equation*}
    \ell_n(f) = -\frac{1}{2}\sum_{i=1}^n ( Y_i - G(f)(X_i))^2
\end{equation*}
and the posterior is then given as
\begin{equation}
\label{eq:generic_posterior}
    \Pi_n (f \mid   \mathcal{D}^{(n)}) \ \propto \exp({\ell_n(f)}) \Pi_n(f).
\end{equation}
For the complete derivation, see \cite{NicklLectureNotes2023} chapter 1.2.3. We can then define a point estimator for $f_0$, given by the map of \cref{eq:generic_posterior}, which is
\begin{equation}
\label{eq:map_of_generic_post}
    f_{\rm MAP} \in \underset{f \in \mathcal{F}}{\rm argmax } \, \Pi_n(f \mid \mathcal{D}^{(n)} ).
\end{equation}
 
The general questions we are interested in are given some prior $\Pi_n(f)$ and data $ \mathcal{D}^{(n)}$, is the posterior $\Pi_n(f \mid  \mathcal{D}^{(n)})$ consistent, and at what rate does it contract to $f_0$, in the sense of \cref{def:contraction}. The second question is, under what conditions does the posterior provide a measure of uncertainty that coincides with the frequentist notion of uncertainty? We remark that from the model \cref{eq:obs}, the resulting posterior will be over $G(f)$; however, we would like to have a posterior of $f$. We can extend the analysis of \cref{eq:generic_posterior} to an induced posterior for $f$, from stability results that come from the forward regularity of the operator $\mathcal{L}_f u_f$, see for example \cite[chapter 2]{NicklLectureNotes2023}.
 
It is typical to choose a Gaussian prior over $f$, see \cite{koers2024, Nickl2020BVM_Schro, NicklWang}, of the form
\begin{equation}
\label{eq:gaussian_prior}
    \Pi_\tau(f) \sim \mathcal{N}(0, \tau_n^{2} \Lambda_\beta^{-1}), \qquad \tau_n^2 := n^{-d/(2\beta + d)},
\end{equation}
where $\Lambda_\beta^{-1} = \rm {diag}(\lambda^\beta_1,...,\lambda^\beta_D)$ and $\beta > 0$ is the smoothness index of $f_0$. The associated penalty $\tau_n^{-2}\|f\|^2_{h^\beta}$ is that of the Tikhonov functional \cref{eq:loss_functional}, of which \cref{eq:map_of_generic_post} is the minimiser, see \cref{subsec:intro_early_stopping}. In this paper the Gaussian prior will instead be placed on the linear parameter $v$, and the prior on $f$ taken to be the induced one, see \cref{subsec:linearisation}.
\begin{remark}
    To achieve the optimal rate, depends on knowing the $\beta$ such that $f_0 \in H^\beta$ \cite{NicklWang} and \cite{NicklLectureNotes2023}. Thus, these results depend on knowing the truth smoothness of $f_0$, something that is not known in practice. Therefore, choosing $\Lambda_\beta$ and the scaling parameter is not possible a priori.
\end{remark}
The question is then whether we can achieve optimal posterior contraction when $\beta$ is unknown.

%%%%%%%%%%%%%%%%%%%%%%%%%%%%%%%%%%%%%%
\subsection{Main Contributions and Outline}
%%%%%%%%%%%%%%%%%%%%%%%%%%%%%%%%%%%%%%
Thus far in the literature, for the non-parametric Bayesian method of PDE-constrained inverse problems, the scale of the prior is either fixed by hand, or choosen by placing a hyperprior. This paper
instead presents a data driven method of early stopping which results in a posterior for $f$ which is shown to contract adaptively and to produce conservative credible sets. Our contributions are as follows:
 
\begin{enumerate}[label=(C\arabic*), leftmargin=*]
\item We extended our analysis to a wider and rather exhaustive class of elliptic semi-linear partial differential equations. 
  We extended the analysis from purely elliptical scenario to parabolic equations in an additional chapter. 
 
\item We extend the method of linearisation beyond Lipschitz solution maps. The linearisation of \cite{koers2024} transfers contraction rates and credible sets from the linearised problem back to the original one only when the solution map $e \colon v \mapsto f$ is Lipschitz. We show in \cref{proposition:transfer} and \cref{proposition:credible} that it suffices for $e$ to admit a modulus of continuity $\omega_n$ (\cref{assumption:modulus}), in which case the rate becomes $\omega_n(\epsilon_n)$ and credible sets have diameter $2\omega_n(2r_n)$. 
 
\item We provide a condition for the existence of such solution map $e$.  We give a sufficient condition directly on the non-linearity of the pde when the nonlinear can be completely written as a function $c(u_f,f)$, see \cref{prop:nemytskii_modulus}. This thus extends our work to the entire class of semilinear elliptic equations.
 
\item We then show that the posterior arising from  the data-driven prior is adaptive up to certain smoothness of the original parameter. Choosing $\tau_n^2 = \tau_{\rm dp}^2$ by the stopping rule \eqref{eq:dp_stopping_time}, the linearised posterior contracts at $\epsilon_n \asymp n^{-2\beta/(2\beta + 2\alpha + 2p + 1 + d)}$, see \cref{lem:linear_post_contracts}, and this rate transfers to the posterior for $f$, see \cref{thrm:orginal_post_contracts}. 
 
\item For fixed smoothness, we show that the credible sets of the resulting posterior for $f$ have asymptotic frequentist coverage one \cref{corollary:coverage_original_post}, with radius
of the same order as in the linearised problem when $c$ is differentiable in $u_f$, and of order $\omega_n(r_{n})$ otherwise. 
 
\item We provide an implementable algorithm and supporting numerics for \cref{eq:schroedinger_eq}. We show that the prior scale parameter can be interpreted as the homotopy parameter in \cite{reich10} identify it
with time in the ensemble Kalman--Bucy filter, so that early stopping becomes a
stopping time for the filter. Our numerical results are for the stationary
Schr\"odinger equation can be found in \cref{subsec:numerical_results} which confirm the
theory of \cref{subsec:results_schro} at finite $n$. The code is  available at \url{https://github.com/Tienstra/EarlyStoppingBIP} .
 
\end{enumerate}
 
\noindent
The paper is organised as follows. In \cref{sec:background_and_prelim} we present 
the background on semilinear elliptic PDEs, the Bayesian model, early stopping
via the discrepancy principle, and the linearisation method of \cite{koers2024}.
\Cref{sec:main_results} presents the general theory, (C1)--(C4). In \cref{sec:examples} we show how the results from \cref{sec:main_results} apply to a variety of examples particularly the stationary and parabolic Allen-Cahn equations, and as edge case stationary Darcy flow.
\Cref{sec:enkbf} presents the ensemble Kalman--Bucy formulation, states the algorithm, and shows the numerical experiments. We conclude in \cref{sec:conclusion}. Auxiliary results used from the literature and additional in the proofs are collected in \cref{sec:appendix_A,sec:appendix_B,sec:appendix_C}.

%%%%%%%%%%%%%%%%%%%%%%%%%%%%%%%%%%%%%%
\subsection{Previous Work}
%%%%%%%%%%%%%%%%%%%%%%%%%%%%%%%%%%%%%%
This work builds on a broad spectrum of existing results, particularly the theory of regularisation in inverse problems, statistical early stopping, Bayesian inverse problems, and empirical prior hyperparameter tuning.

The regularisation of inverse problems, especially Tikhonov-type regularisation with hyperparameter selection via the discrepancy principle, has been thoroughly studied in the literature. See \cite{EngHanNeu96} for a comprehensive treatment of the linear inverse problem setting with bounded noise. In a related direction, \cite{Hanke1995} analyses early stopping for gradient descent using a discrepancy-based stopping rule in the nonlinear setting. Their study focuses on mildly ill-posed deterministic inverse problems, requiring the initialisation to be sufficiently close to the ground truth to ensure local convexity of the Tikhonov functional.

Statistical early stopping also has a rich body of literature. For instance, \cite{BlanchardHoffmannReiss-bis, Stankewitz} investigated early stopping strategies for statistical linear inverse problems using truncated SVD. Further, \cite{BlanchardHoffmannReiss} extend this to discrepancy-based stopping rules for both gradient descent and Tikhonov regularisation in linear settings. The recent work \cite{tienstra2025} generalises these results by incorporating regularisation operators into the penalisation term, providing a Bayesian interpretation of the stopping rule.

Bayesian inverse problems, both in linear \cite{knapik2011, Stuart} and nonlinear settings \cite{Nickl2020BVM_Schro, NicklLectureNotes2023, giordano2020consistency, NicklWang}, have also been extensively developed. In particular, hyperparameter selection for Gaussian priors in linear Bayesian inverse problems has been approached both empirically and hierarchically in \cite{szabo2013empirical}.

Finally, the work of \cite{koers2024}, which provides a framework to reparameterize nonlinear inverse problems into linear inverse problems, enables the direct application of the theoretical results from \cite{tienstra2025}. However, this general method requires a case-by-case checking of conditions. This linearisation enables the application of exact methods such as the Ensemble Kalman Filter (EnKF), which can evolve the prior distribution dynamically toward the true posterior. Building on the homotopy approach formulated in \cite{reich10}, the scale parameter of the prior covariance can be interpreted as a time-like parameter, allowing for a continuous deformation of the prior into the posterior, thus providing a Bayesian iterative method to compute the target posterior distribution. 
%%%%%%%%%%%%%%%%%%%%%%%%%%%%%%%%%%%%%%
\subsection{Notation}
%%%%%%%%%%%%%%%%%%%%%%%%%%%%%%%%%%%%%%
We further denote $L^2(\mathcal{O})$ as the space of square integrable functions on $\mathcal{O}$,  with inner product $\langle \cdot , \cdot \rangle_{L^2(\mathcal{O})}$ induced norm $\|\cdot \|_{{L^2(\mathcal{O})}}$. We further denote by $h^\beta$ the sequence space
\begin{equation}
\label{eq:h_beta}
    h^{\beta} := \{f \in \ell^2(\mathbb{N}): \|f\|^2_{h^\beta} = \sum_i^\infty \lambda_i^\beta f_i^2 < \infty\}.
\end{equation}
For $s>0$ and $\mathcal{O}\subset\mathbb{R}^d$ we further define the (dimension-adjusted) sequence-space Sobolev space
\begin{equation}
\label{eq:H_beta_d}
    H_d^s := \Big\{ v = (v_i)_{i\in\mathbb{N}} \in \ell^2(\mathbb{N}) : \|v\|_{H_d^s}^2 := \sum_{i=1}^\infty i^{2s/d}\, v_i^2 < \infty \Big\},
\end{equation}
so that $\beta'$ in $v_0 \in H_d^{\beta'}$ denotes the (true, unknown) regularity of the ground truth $v_0$, and, writing $\widetilde{v}_0 := C_0^{-1/2}v_0$ for the whitened parameter, we define $\beta$ to be such that $\widetilde{v}_0 \in H_d^\beta$; see \cref{assumption:1}.

We define the following additional standard statistical notation, see \cite{Ghosal2017book}. For two numbers $a$ and $b$, we denote the minimum of $a$ and $b$ by $a \wedge b$. For two sequences $(a_n)_n$ and $(b_n)_n$ in $\mathbb{R}_{+}$, $a_n \lesssim b_n$, respectively $a_n \gtrsim b_n$ denote inequalities up to a multiplicative constant. $a_n \asymp b_n$ denotes that $a_n \lesssim b_n$ and $a_n \gtrsim b_n$ hold. $\ell^2(\mathbb{N})$ denotes the space of sequences that are square summable with index $i \in \mathbb{N}$, and its norm is denoted by $\|a \|_{\ell^2(\mathbb{N})} = \left(\sum_i a_i^2\right)^{1/2}$. Unless stated otherwise, $d$ denotes a semi-metric on $\mathcal{F}$; in our examples this will be $d(f,f') := \|f-f'\|_{L^2(\mathcal{O})}$. Finally when we write
\begin{equation*}
    \Pi_n(\mathcal{B}_n \mid Y_n) \overset{P_{f_0}}{\rightarrow} 1
\end{equation*}
for the set $B_n = \{x \in \mathcal{F} \mid d(x,f_0) \leq \epsilon_n \}$, observations $Y_n$, and $P_{f_0}$ the law of the data $Y_n$ under the true, data-generating parameter $f_0$, we mean that
\begin{equation}
 P_{f_0}\left( \Pi(x : d(x,f_0) \geq  \epsilon_n \mid Y_n) > \delta_n \right) \rightarrow 0
\end{equation}
as $n \rightarrow \infty$ for every $\epsilon_n,
\delta_n \rightarrow 0$. That is the posterior concentrates around the ball that shrinks to the truth.

%%%%%%%%%%%%%%%%%%%%%%%%%%%%%%%%%%%%%%
\section{Background Theory}\label{sec:background_and_prelim}
%%%%%%%%%%%%%%%%%%%%%%%%%%%%%%%%%%%%%%
In this section we review important background theory for \cref{sec:main_results}. We first give a short introduction to semilinear pdes, and then introduce early stopping in the Bayesian context and a linearisation method for non-linear pdes. 
\subsection{Semilinear elliptic partial differential equations}
Let $\mathcal{O}$ be a bounded domain in $\mathbb{R}^d$ with a $C^1$ boundary.
Let $U\subset H^{2m}(\mathcal{O})$ be a $H^{2m}$ closed subspace and $V = L^2(\mathcal{O})$.

Let the parameter space $\mathcal{F}$ be a closed subspace of $H^{l}(\mathcal{O})$.
We then consider differential operators of the form
\begin{equation}
\mathcal{L}_f(u) = \mathbb{L} u - C(u; f)
\end{equation}
where $\mathbb{L}:U \rightarrow V$ is a symmetric $U$-coercive differential operator of order $2m$
and the Nemytskii operator is given pointwise by
\begin{equation}
C(u; f)(x) = c(x, J^k u(x); J^l f(x)).
\end{equation}
 $J^k u := (\partial^\alpha u)_{|\alpha| \leq k} \in \mathbb{R}^M$ and $J^l f$ are the $k$-th order jet prolongation of $u$ and $l$-th order jet prolongation of f, respectively, where
$\alpha$ is a multi-index of degree $|\alpha| \leq k < 2m$.

Based on the Sobolev embedding $H^{s-|\alpha|}(\mathcal{O}) \hookrightarrow L^\infty(\mathcal{O})$, for any jet depth $j \ge 0$ and Sobolev regularity $s > 0$, we define the local and global multi-index sets:
\begin{align*}
    A_{\mathrm{loc}}(j, s) &= \{ \alpha : |\alpha| \leq j \text{ and } s - |\alpha| > d/2 \}, \\
    A_{\mathrm{glob}}(j, s) &= \{ \alpha : |\alpha| \leq j \text{ and } s - |\alpha| \leq d/2 \}.
\end{align*}
We apply this partition to the state jet $J^k u$ using $s=2m$, and to the parameter jet $J^l f$ using $s=l$. 
For vectors $\xi \in \mathbb{R}^M$ and $\zeta \in \mathbb{R}^{M_f}$, their restrictions to the local indices are denoted by $\xi_{\mathrm{loc}} := (\xi_\alpha)_{\alpha \in A_{\mathrm{loc}}(k, 2m)}$ and $\zeta_{\mathrm{loc}} := (\zeta_\alpha)_{\alpha \in A_{\mathrm{loc}}(l, l)}$.

To cover a general class of semilinear PDEs, we adapt standard growth conditions from nonlinear functional analysis (see, e.g. \cite{zeidler_nonlinear_1990}). 
By explicitly partitioning the jet indices into those that embed continuously into $L^\infty$ ($A_{\mathrm{loc}}$) and those that only embed into Lebesgue spaces ($A_{\mathrm{glob}}$), we formulate a unified set of conditions on the Nemytskii operator. 
This structure allows us to rigorously establish the Lipschitz continuity of the inverse solution map in the subsequent lemmas.

\begin{assumption}
\label{ass:nonlinearity}
Let $c:\mathcal{O}\times\mathbb{R}^{M}\times\mathbb{R}^{M_f}\to\mathbb{R}$ be a Carath\'eodory function. 
\begin{enumerate}[label=(C\arabic*), leftmargin=*]

\item For every $R,R_f>0$, there exist $a_{R,R_f}\in L^2(\mathcal{O}), b_{R,R_f} \geq 0$, and exponents $q_\alpha \geq 1$ and $p_\alpha \geq 1$, such that
\begin{equation}
\label{eq:growth_bound}
|c(x,\xi;\zeta)| \leq a_{R,R_f}(x) + b_{R,R_f} \left( \sum_{\alpha\in A_{\mathrm{glob}}(k, 2m)} |\xi_\alpha|^{q_\alpha} + \sum_{\alpha\in A_{\mathrm{glob}}(l, l)} |\zeta_\alpha|^{p_\alpha} \right)
\end{equation}
for almost every $x\in\mathcal{O}$, for every $\xi\in \mathbb{R}^M$ with $ |\xi_{\mathrm{loc}}|\leq R$, and every $\zeta\in \mathbb{R}^{M_f}$ with $|\zeta_{\mathrm{loc}}|\leq R_f$.
Furthermore, we assume $q_\alpha \leq \frac{d}{d-2(2m-|\alpha|)}$ and $p_\alpha \leq \frac{d}{d-2(l-|\alpha|)}$.

\item For every $R,R_f>0$, there exist a non-decreasing function $ \omega_{R,R_f}:[0,\infty)\to[0,\infty)$ with $\omega_{R,R_f}(0)=0$ where $t \mapsto \omega_{R,R_f}(\sqrt{t})^2$ is concave, a constant $L_{R,R_f}<\infty$, and exponents $ \gamma_\alpha\in(0,1]$ and $\delta_\alpha\in(0,1]$, such that
\begin{equation}
\label{eq:continuity_bound}
\begin{aligned}
|c(x,\xi_1;\zeta_1)-c(x,\xi_2;\zeta_2)| \leq & \, L_{R,R_f} \omega_{R,R_f} (|\xi_{1,\mathrm{loc}}-\xi_{2,\mathrm{loc}}| + |\zeta_{1,\mathrm{loc}}-\zeta_{2,\mathrm{loc}}|)\\
&+ L_{R,R_f} \sum_{\alpha\in A_{\mathrm{glob}}(k, 2m)} \left( |\xi_{1,\alpha}|^{q_\alpha-\gamma_\alpha} + |\xi_{2,\alpha}|^{q_\alpha-\gamma_\alpha} \right) |\xi_{1,\alpha}-\xi_{2,\alpha}|^{\gamma_\alpha}\\
&+ L_{R,R_f} \sum_{\alpha\in A_{\mathrm{glob}}(l, l)} \left( |\zeta_{1,\alpha}|^{p_\alpha-\delta_\alpha} + |\zeta_{2,\alpha}|^{p_\alpha-\delta_\alpha} \right) |\zeta_{1,\alpha}-\zeta_{2,\alpha}|^{\delta_\alpha}
\end{aligned}
\end{equation}
for almost every $x\in\mathcal{O}$, and for all $\xi_i,\zeta_i$ subject to the local bounds $R, R_f$.

\item For a given $f_0\in\mathcal F$ and corresponding solution $u_{f_0}\in U$, there exist $r,s>0$, an invertible bounded linear operator $A\in L(\mathcal F,V)$ with $\|A^{-1}\|< \infty$, and $\theta\in(0,1)$ such that
\begin{equation}
\label{eq:param_approx}
\left\| C(u;f_1) - C(u;f_2) - A(f_1-f_2) \right\|_V \leq \frac{\theta}{\|A^{-1}\|} \|f_1-f_2\|_{\mathcal F}
\end{equation}
for all $u\in B_r(u_{f_0})$ and $f_1,f_2\in B_s(f_0)$.
\end{enumerate}
\end{assumption}

\begin{proposition}
\label{prop:nemytskii_modulus}
Suppose \cref{ass:nonlinearity} holds true,
then the Nemytskii operator $C:U\times\mathcal F\to V$
is well defined and maps bounded subsets of $U\times\mathcal F$ into bounded subsets of $V$.
Moreover, let $f_0\in\mathcal F$ and $u_{f_0}\in U$ be fixed. Then there exist constants $r,s>0$ and a modulus of continuity $\omega_C:[0,\infty)\to[0,\infty)$ with $\omega_C(0)=0$ 
such that
\begin{equation}
\label{eq:C_modulus}
\| C(u;f)-C(v;f)\|_{V} \leq \omega_C(\|u-v\|_U)
\end{equation}
for all $ u,v\in B_r(u_{f_0})$,$ f\in B_s(f_0)$.
More specifically, there exist constants $w_0,w_1>0$ satisfying,
\begin{equation}
\label{eq:explicit_C_modulus}
\omega_C(t) = w_0 \left(\omega_{R,R_f}(w_1t) + \sum_{\alpha\in A_{\mathrm{glob}}} t^{\gamma_\alpha} \right).
\end{equation}
Proof in \cref{sec:appendix_C}.
\end{proposition}
\begin{remark}\label{rem:frechet_derivative}
  Assume that the induced Nemytskii operator $(u, f) \mapsto C(u; f)$ admits a partial Fréchet derivative $\partial_f C(u; f) \in L(\mathcal{F}, V)$ for all $(u, f) \in B_r(u_{f_0}) \times B_s(f_0)$. 
  If the linear operator $\partial_f C(u_{f_0}; f_0)$ has a bounded inverse, and the mapping $(u, f) \mapsto \partial_f C(u;f)$ from $U \times \mathcal{F}$ into $L(\mathcal{F}, V)$ is continuous at $(u_{f_0}, f_0)$, then condition (C3) is satisfied.
  Proof in \cref{sec:appendix_C}.
\end{remark}
%\begin{remark}\label{rem:unified_bounds}
%  Because $U \subset H^{2m}(\mathcal{O})$, the components $D^\alpha u$ for $\alpha \in A_{\mathrm{loc}}$ compactly embed into $L^\infty(\mathcal{O})$, ensuring their pointwise inputs are globally bounded by some $R > 0$. 
%  If $2m - k > d/2$ (e.g., $d \leq 3$ for $k=0$), then $A_{\mathrm{glob}}$ is empty. 
%  The polynomial sum in \cref{eq:growth_bound} vanishes, recovering purely local $L^2$-boundedness for fast-growing physical non-linearities. 
%  Conversely, if $A_{\mathrm{glob}}$ is populated, the polynomial condition is rigorously enforced only on the components that lack the $L^\infty$ embedding.
%\end{remark}
%
%\begin{remark}\label{rem:derivative_loss}
%  If the parameter $f$ is differentiated in the PDE ($l \geq 1$), the Fréchet derivative $A$ becomes a differential operator. By equipping $\mathcal{F} \subseteq W^{l,\infty}(\mathcal{O})$ with a sufficiently strong topology, we mathematically absorb the associated loss of derivatives. This ensures the inverse operator $A^{-1}: V \to \mathcal{F}$ remains bounded, allowing condition (C3) to hold.
%\end{remark}
Henceforth, we denote the compact self adjoint inverse to $\mathbb{L}$ \cite[Chapter 6]{evans_partial_2010} by $\mathcal{K} := \mathbb{L}^{-1}:V\rightarrow U$.
\begin{remark}
\label{rem:spillting_pde}
  For a differential operator of order $2m$, the normal trace operators $\frac{\partial^j u}{\partial \nu^j}$ mapping $H^{2m}(\mathcal{O})\to H^{2m-j-\frac{1}{2}}(\partial \mathcal{O})$ for $j = 0, \dots, m-1$ are well defined and bounded. 
  Without loss of generality, we restrict ourselves to the case of homogeneous Dirichlet boundary conditions, i.e.
\begin{equation}
\label{eq:generic_pde}
    \left\{ 
    \begin{alignedat}{4}
        \mathbb{L}\tilde{u}_f - C(\tilde{u}_f + \tilde{g}; f)& = & \,  h &  \quad \text{on } \mathcal{O},\\
        \frac{\partial^j \tilde{u}_f}{\partial \nu^j} &=& 0 &  \quad \text{on } \partial\mathcal{O} \text{ for } j=0,\dots,m-1.
    \end{alignedat}
    \right.
\end{equation}
where $\tilde{u}_f = u_f - \tilde{g}$, and $\tilde{g}$ is the unique solution to the linear boundary value problem
\begin{equation}
    \label{eq:generic_pde_boundary}
    \left\{ 
    \begin{alignedat}{4}
        \mathbb{L}\tilde{g} &=& 0 &  \quad \text{on } \mathcal{O},\\
        \frac{\partial^j \tilde{g}}{\partial \nu^j} &=& g_j &  \quad \text{on } \partial\mathcal{O} \text{ for } j=0,\dots,m-1.
    \end{alignedat}
    \right.
\end{equation}%
By absorbing the boundary data into the operator as a fixed shift $\tilde{g}$, we can exclusively study the homogeneous problem. 
This observation allows us to choose $U = H_{0}^{m}(\mathcal{O})\cap H^{2m}(\mathcal{O})$, which is a closed subspace of $H^{2m}(\mathcal{O})$.
\end{remark}

\subsection{Linearisation of  Non-linear Inverse Problems}\label{subsec:linearisation}
An integral part of extending the results of \cite{tienstra2025} is the linearization method found in \cite{koers2024}. We give an overview of their approach below and point the reader to the source for the remaining details. The approach involves using the splitting in \cref{eq:generic_pde} to define an inverse in which we can arrive at a linear problem. From now on, let $h=0$ in \cref{eq:generic_pde}. As $h$ is a known function in $V$ that is additive we can let $c(u,f) = c(u,f) -h$. We now consider the continuous  observations  
\begin{equation}
\label{eq:cont_obvs}
    Y_{n} = u_f + n^{-1/2} \Xi
\end{equation}
where this should now be understood as a process similar to what was defined in \cref{subsec:intro_early_stopping}. Recall \cref{eq:generic_pde}, and that $\mathcal{K}$ inverse of $\mathbb{L}$ and $\tilde{g}$ was such that \cref{eq:generic_pde_boundary} holds. 
As $\tilde{g}$ is unique, we can write the solution of \cref{eq:generic_pde_full_operator} as
\begin{equation}
    u_f = \mathcal{K}\mathbb{L}u_f + \tilde{g}
\end{equation}
respectively, it holds that
\begin{equation}
    \mathbb{L}(\mathcal{K} \mathbb{L}u_f + \tilde{g}) = \mathbb{L}\mathcal{K}(\mathbb{L}u_f) + 0 = \mathbb{L} u_f
\end{equation}
on $\mathcal{O}$. 
Let $v=\mathbb{L}u_f$ we can then define continuous observations 
\begin{align}
\label{eq:linear_model}
    \widetilde{Y}_n &:= Y_n - \tilde{g} = \mathcal{K}(\mathbb{L}u_f) + \frac{1}{\sqrt{n}} \Xi \nonumber \\
    &:= \mathcal{K}v + \frac{1}{\sqrt{n}} \Xi
\end{align}
\begin{equation}
    \label{eq:solution_map_e}
    f = e(\mathbb{L} u_f).
\end{equation}
where, following the notations of \cite{koers2024}, we now consider two different posterior distributions; the posterior arising from the non-linear problem \cref{eq:generic_pde_full_operator} with $h=0$ which is given and denoted as
\begin{equation}
    \label{eq:original_post}
    \Pi_n (f \mid Y_{n}) \ \propto \  \Pi_n(f) L(Y_n \mid u_f)
\end{equation}
where $\Pi_n(f)$ denotes the prior of $f$, and $L(Y_{n} \mid u_f)$ denotes the likelihood of $Y_{n} \mid u_f$ under what we refer to as model (N)
\begin{equation}
\label{eq:nonlinear_model}
    Y_{n} = u_f + n^{-1/2} \Xi
\end{equation}
Similarly, the posterior arising from the linear problem, which is given and denoted as
\begin{equation}
    \label{eq:linearized_post}
    \widetilde{\Pi}_n (v\in \cdot \mid \tilde{Y_{n}}) \ \propto \  \widetilde{\Pi}_n(v) L(\widetilde{Y}_{n} \mid \mathcal{K}v)
\end{equation}
where $\widetilde{\Pi}_n(v)$ denotes the prior of $v$, and $L(\widetilde{Y}_{n} \mid \mathcal{K}v)$ denotes the likelihood of $\widetilde{Y}_{n} \mid \mathcal{K}v$. 
\begin{remark}
    In \cite{koers2024}, the goal is to do the frequentist analysis for the Gaussian posterior \cref{eq:linearized_post}, and pull back the results to the \textit{original} posterior \cref{eq:original_post}. We remark that this \textit{original} posterior is that which arises from the induced prior $\Pi(v)$. That is
    \begin{equation}
        \Pi_n(f) := \widetilde{\Pi}_n(e(v))
    \end{equation}
    \begin{align}
        \widetilde{\Pi}_n(e(v) \mid \widetilde{Y}_n) \ &\propto \  \widetilde{\Pi}_n(e(v))L(\widetilde{Y}_n \mid v) \\
        & \propto \ \Pi_n(f)L(\widetilde{Y}_n \mid v) \\
        & \propto \ \Pi_n(f)L(Y_n \mid f) 
    \end{align}
    where $e$ is the solution map $f \mapsto \mathbb{L}u_f$ and $v = \mathbb{L}u_f$. The posterior distribution $\widetilde{\Pi}_n(e(v)$ is the induced posterior from the linear Gaussian one \cref{eq:linearized_post} via the map $e$. In this paper, we consider the original posterior to be such that the above equations hold. 
\end{remark}
The key questions in frequentist Bayesian analysis: under what conditions does the posterior contract to the ground truth function, at what rate, and under what conditions does the posterior spread coincide with frequentist confidence intervals, can then be easily asked for linear posterior \cref{eq:linearized_post}.  The key result of \cite{koers2024} is that asymptotically, the above-mentioned theoretical results of \cref{eq:linearized_post} can be pulled back to the original posterior  \cref{eq:original_post}. It seems plausible then that if we choose a prior such for $v$ that is $\mathcal{N}(0, \tau_n^2 C_0)$, that we could choose $\tau_n^2$ via early stopping, see  \cref{subsec:intro_early_stopping}, and pull back results of the posteior which would now depend on $\tau_{\rm dp}$
\begin{equation}
    \widetilde{\Pi}_{n, \tau_{\rm dp}} (v \mid \widetilde{Y}_n)
\end{equation}
to the original posterior \cref{eq:original_post}.

\subsection{Early Stopping in the Bayesian context}\label{subsec:intro_early_stopping}

Early stopping is generally applied in the non-Bayesian setting, for example, to Tikhonov regularisation, see \cite{EngHanNeu96, BlanchardHoffmannReiss}. However, it is known that Tikhonov regularisation and the Bayesian setting are intrinsically linked \cite{Stuart}, via computing the MAP estimator, which is the minimiser of the Tikhonov functional. More specifically, we have that the MAP estimator of the posterior, denoted by $v_{\rm map}$, is the minimiser of
\begin{equation}
 \label{eq:loss_functional}
     \mathcal{T}_{(\tau_n)}(v) =   \|P_n(\mathcal{K}v - \widetilde{Y}_n)\|^2_2 +
\tau_{n}^{-2} ||C_0^{-1/2} v||^2_1 .
\end{equation}
where $P_n$ is an appropriate projection operator onto a $D(n)$-dimensional subspace of $V$ (see \cref{eq:projector} below), and  $C_0$ is the covariance operator of the prior and now $\mathcal{K}$ is a linear forward operator in place of $G$ in \cref{eq:obs}. For an estimator $\widehat{v}$, the residuals are defined as
\begin{equation}
    \label{eq:residuals_projected}
     R_{\tau_{n}} := || P_n(\widetilde{Y}_n - \mathcal{K} \widehat{v}_{\tau_{n}})||^2.
\end{equation}
Suppose further that we have an iterative method such that for each $\tau_n \in \mathbb{R}_{+} \cup \{0\}$ we can construct a sequence of estimators
\begin{align*}
    (\widehat{v}_{\tau_n})_{\tau_n}
\end{align*}
such that they minimise \cref{eq:loss_functional},
\begin{align*}
    \widehat{v}_{\tau_n} = \text{argmin } \mathcal{T}_{(\tau_n)}(v).
\end{align*}
and can be ordered in decreasing bias and increasing variance. Suppose also that we choose $P_n$ such that $\widetilde{Y}_n$ is projected to be of dimension $D(n) \leq n$. Then for each $\tau_n \mapsto \widetilde{v}_{\tau_n}$ we can stop the iterative process at
\begin{align}
\label{eq:dp_stopping_time}
    \tau_{\rm dp}(n) := \inf\,\{\tau_{n} > 0 : R_{\tau_{n}} \leq \kappa \}, \text{ where } \kappa \asymp D(n)/n.
\end{align}
When the noise level is known and constant, \cref{eq:dp_stopping_time} is called the discrepancy principle, see \cite{EngHanNeu96}. In \cite{tienstra2025} it was shown, in the bayesain linear setting, that we can choose the optimal scaling parameter, $\tau_n$, of \cref{eq:loss_functional} according to the stopping rule \cref{eq:dp_stopping_time} for appropriately chosen $C_0$ and for $\kappa \asymp D(n)/n$ where $D(n)$ is derived below. To do this, one must project $Y_i$ into some finite-dimensional subspace to define the stopping criterion. From \cite{tienstra2025}, we know that the appropriate $D(n)$ depends on the effective dimension and should be chosen as
\begin{equation}
    \label{eq:project_dim}
    D(n) \asymp n^{1/(2p+d)}.
\end{equation}
where $p$ is the decay of the eigenvalues of the forward operator and one assumes that with respect to some basis $(t_i)_{i\in\mathbb{N}}$ the eigenvalues $(\kappa_i)_{i \in \mathbb{N}}$ of the forward operator are such that
 \begin{equation}
     \label{eq:decay_of_kappa}
     \kappa_i \asymp i^{-p/d} \quad i \in \mathbb{N}.
 \end{equation}
We then observe \cite{tienstra2025}
\begin{equation}
    \label{eq:projector}
    \langle P_n \widetilde{Y}_n,t_i \rangle_2 = \left\{ \begin{array}{ll}
    \widetilde{Y}_{n,i} & {\rm if}\,\,i\le D(n)\\
    0 & {\rm otherwise}
    \end{array} \right.
\end{equation}
And our observations in sequence space are then
\begin{equation}
    \label{eq:truncated_model}
    \widetilde{Y}_i = \kappa_i v_{0,i} + n^{-1/2}\xi_i \quad i=1,...,D(n) \ \forall n \in \mathbb{N}
\end{equation}
where $(\xi_i)_{i\in\mathbb{N}}$ are i.i.d.\ $\mathcal{N}(0,1)$, the coordinates of the Gaussian white noise $\Xi$ of \cref{eq:linear_model} with respect to the basis $(t_i)$.
Then by Theorem 2.1 in \cite{tienstra2025} (listed in appendix as \cref{thrm:tien_contraction_rate_early_stopping} ), the posterior $\widetilde{\Pi}_n(v \mid \widetilde{Y}_n)$ is such that for $\epsilon_n \asymp n^{-2\beta/(2\beta + 2p + 2\alpha + 1 + d)}$, and $M_n \rightarrow \infty$
\begin{equation}
     \widetilde{\Pi}_{n, \tau_{\rm dp}} (v\in V : d_V(v, v_0) \geq M_n \epsilon_n \mid \widetilde{Y}_n) \overset{P^n_{v_0}}{\longrightarrow} 0
\end{equation}
holds. Furthermore it was shown, in corollary 2.1 of \cite{tienstra2025}, see \cref{corollary:tien_assympotic_coverage} in appendix for reference, that the coverage, see \cref{def:credible_set}, tends to 1. That is
\begin{equation}
    \widetilde{\Pi}_{n, \tau_{\rm dp}}(v_0 \in v_{\rm map} + B(r_{n,c}) \mid \widetilde{Y}_n) \rightarrow 1.
\end{equation}
as $n \rightarrow \infty$. We would now like to extend such results to posteriors \cref{eq:generic_posterior}, which arise from \cref{eq:generic_pde}. To this, we use a linearisation scheme, which is the topic of the next section.

%%%%%%%%%%%%%%%%%%%%%%%%%%%%%%%%%%%%%%

\section{General Theoretical Results}\label{sec:main_results}

%%%%%%%%%%%%%%%%%%%%%%%%%%%%%%%%%%%%%%
This section first extends the linearisation method discussed in \cref{subsec:linearisation} to more general solution maps $e$, see \cref{def:mod_cont}, particularly when $e$ is $\rho-$ Hölder continuous with $\rho < 1$ see \cref{proposition:transfer} and \cref{proposition:credible}. We then present a condition on $c$ for which we are guaranteed, that the solution map $e$, admits a modulus of continuity, and in which case we are guaranteed that it will be Lipschitz on certain open sets, see \cref{lem:c_continuous} and \cref{prop:nemytskii_modulus} respectively. We then show in \cref{thrm:orginal_post_contracts} and \cref{corollary:coverage_original_post}, that if the scaling parameter of \cref{eq:linearized_post} is chosen via early stopping \cref{eq:dp_stopping_time}, the posterior \cref{eq:original_post} arising from \cref{eq:nonlinear_model} contracts optimally for the reparametrized problem \cref{assumption:1}, and has asymptotic frequentist coverage equal to 1.

\begin{definition}
[Modulus of continuity]
\label{def:mod_cont}
A modulus of continuity is a function $\omega : [0, \infty] \rightarrow [0,\infty)$  such that $\omega(0)=0$, and is used to measure quantitatively the uniform continuity of functions. We say that $f$ admits $\omega$ as a modulus of continuity if

\begin{equation}
|f(x) - f(y)| \le \omega(|x-y|)
\label{eq:mod_cont}
\end{equation}

for all $x,y$ in the domain of $f$.
\end{definition}
\begin{assumption}[Modulus of continuity of the solution map]
\label{assumption:modulus}
There exist Borel sets $V_n \subseteq V$ and functions
$\omega_n \colon [0,\infty) \to [0,\infty)$ such that:
\begin{enumerate}
  \item[(i)] each $\omega_n$ is non-decreasing with
    $\omega_n(t) \to 0$ as $t \downarrow 0$;
  \item[(ii)] for every $v \in V_n$,
    \begin{equation}\label{eq:modulus}
      \bigl\| e(v) - e(v_0) \bigr\|_{L^2}
        \;\leq\; \omega_n\bigl( \|v - v_0\| \bigr);
    \end{equation}
  \item[(iii)] $\widetilde{\Pi}_n\bigl( V_n \mid \widetilde{Y}_n \bigr)
    \xrightarrow{P_{v_0}} 1$.
\end{enumerate}
In our constructions $V_n$ is taken to be a ball around the truth (e.g.\ $V_n = B^V_n(v_{f_0})$ as in \cref{eq:ball_of_lip_solution_map}), so $v_0 \in V_n$ for every $n$; condition (iii) is the additional requirement that the posterior concentrate on this set.
If in addition each
$\omega_n$ is concave we have that the $t \mapsto \omega_n(t)/t$ is
non-increasing, and so then
\begin{equation}\label{eq:concave}
  \omega_n(ct) \;\leq\; c \, \omega_n(t),
  \qquad c \geq 1,\ t \geq 0 .
\end{equation}
\end{assumption}
 
 \begin{definition}
 \label{def:contraction}
    Let  $(\epsilon_n)_n$ be a sequence of positive numbers. Then $(\epsilon_n)_n$ is a posterior contraction rate at the  parameter $G(f_0)$ wrt to the some semi-metric $d$ if for every sequence $(M_n)_n \rightarrow \infty$, it holds that,
    \begin{align}
        \label{eq:contraction_rate}
        \Pi_n (f\in \mathcal{F} : d(f, f_0) \geq M_n \epsilon_n \mid \mathcal{D}^{(n)} ) \overset{P^n_{f_0}}{\longrightarrow} 0
    \end{align}
    as $n \rightarrow \infty$. Where $\Pi_n(\cdot \mid \mathcal{D}^{(n)} )$ is the posterior given observations $Y$ and given prior $\Pi_n$. The maximum such $\epsilon_n$ that \cref{eq:contraction_rate} holds is called the posterior contraction rate. If $\epsilon_n$ matches the minimax rate of $G(f_0)$, then the posterior contracts optimally to $G(f_0)$.
\end{definition}

\begin{proposition}[Pull back contraction rate]
\label{proposition:transfer}
Suppose the posterior distribution $\widetilde{\Pi}_n(v \in \cdot \mid
\widetilde{Y}_n)$ of $v$ in model (L) contracts under $v_0$ to
$v_0 = \mathbb{L}u_{f_0}$ at rate $\epsilon_n$ in $(V, \|\cdot\|)$.
Let $f = e(v)$ hold and suppose
Assumption~\ref{assumption:modulus} is satisfied. Then, for every sequence
$M_n \to \infty$,
\begin{equation}\label{eq:conclusion}
  \Pi_n\Bigl( f \colon \|f - f_0\|_{L^2}
    > \omega_n(M_n \epsilon_n) \,\Bigm|\, Y_n \Bigr)
  \;\xrightarrow{P_{f_0}}\; 0 .
\end{equation}
If in addition $(\omega_n)$ is concave, then the posterior
distribution of $f$ in model (N) attains the contraction rate
$\omega_n(\epsilon_n)$ under $f_0$ relative to the $L^2$-norm.
\end{proposition}
 
\begin{proof}
This proof follows closely the proof of \cite[proposition. 2.1]{koers2024}
By assumption the prior distributions of $\mathbb{L}u_f$ and of $v$
coincide, and by \eqref{eq:linear_model} and \eqref{eq:solution_map_e} the
conditional distributions of $Y_n - \tilde{g} = \mathcal{K}\mathbb{L}u_f
+ n^{-1/2}\Xi$ given $\mathbb{L}u_f$ and of $\widetilde{Y}_n =
\mathcal{K}v + n^{-1/2}\Xi$ given $v$ coincide as well, under the
identification $v = \mathbb{L}u_f$. Consequently the conditional
distributions of $\mathbb{L}u_f$ given $Y_n - \tilde{g}$ and of $v$
given $\widetilde{Y}_n$ coincide: if $(y, B) \mapsto L(y, B)$ is a
Markov kernel with $v \mid \widetilde{Y}_n \sim L(\widetilde{Y}_n,
\cdot)$, then $\mathbb{L}u_f \mid Y_n \sim L(Y_n - \tilde{g}, \cdot)$.

Let $M_n$ be such that $M_n \to \infty$ and and define 
\[
  \rho_n := \omega_n(M_n \epsilon_n),
  \qquad
  A_n := \bigl\{ v \in V_n \colon
    \|e(v) - e(v_0)\|_{L^2} > \rho_n \bigr\} .
\]
By the identity
$f = e(v)$, we have that $f - f_0 = e(\mathbb{L}u_f) - e(v_0)$, and so
\[
  \Pi_n\bigl( f \colon \|f-f_0\|_{L^2} > \rho_n,\;
    \mathbb{L}u_f \in V_n \bigm| Y_n \bigr)
  \;=\; L\bigl( Y_n - \tilde{g},\, A_n \bigr) ,
\]
and, because $Y_n - \tilde{g}$ is distributed as $\widetilde{Y}_n$
under $v_0$, the right-hand side is distributed as
$\widetilde{\Pi}_n( A_n \mid \widetilde{Y}_n )$. Note that $A_n$ is
defined as a preimage under $e$, so that no injectivity of $e$ is
required at this step.
 
\begin{equation}\label{eq:inclusion}
  A_n \;\subseteq\; \bigl\{ v \colon \|v - v_0\| \geq M_n \epsilon_n
  \bigr\} .
\end{equation}
We prove this claim of \cref{eq:inclusion} by  contradiction. Let $v \in A_n$. Suppose that
$\|v - v_0\| < M_n \epsilon_n$. Since $v \in V_n$, inequality
\eqref{eq:modulus} together with the monotonicity of $\omega_n$ gives
\[
  \|e(v) - e(v_0)\|_{L^2}
  \;\leq\; \omega_n\bigl( \|v - v_0\| \bigr)
  \;\leq\; \omega_n( M_n \epsilon_n )
  \;=\; \rho_n ,
\]

We then have that 
\begin{align*}
  \Pi_n\bigl( \|f - f_0\|_{L^2} > \rho_n \bigm| Y_n \bigr)
  &\leq \Pi_n\bigl( \|f-f_0\|_{L^2} > \rho_n,\;
     \mathbb{L}u_f \in V_n \bigm| Y_n \bigr)
     + \Pi_n\bigl( \mathbb{L}u_f \notin V_n \bigm| Y_n \bigr) \\
  &= \widetilde{\Pi}_n\bigl( A_n \bigm| \widetilde{Y}_n \bigr)
     + \widetilde{\Pi}_n\bigl( V_n^{\,c} \bigm| \widetilde{Y}_n \bigr) \\
  &\leq \widetilde{\Pi}_n\bigl( v \colon \|v - v_0\|
     \geq M_n \epsilon_n \bigm| \widetilde{Y}_n \bigr)
     + \widetilde{\Pi}_n\bigl( V_n^{\,c} \bigm| \widetilde{Y}_n \bigr),
\end{align*}
where the two terms on the right tend to zero in $P_{v_0}$-probability
by the assumed contraction rate and by
Assumption~\ref{assumption:modulus}(iii) respectively. 
 
Finally, suppose $(\omega_n)$ is concave and let $M_n \to \infty$;
without loss of generality $M_n \geq 1$. By \eqref{eq:concave},
$\omega_n(M_n \epsilon_n) \leq M_n \omega_n(\epsilon_n)$, whence
\[
  \bigl\{ \|f - f_0\|_{L^2} > M_n \, \omega_n(\epsilon_n) \bigr\}
  \;\subseteq\;
  \bigl\{ \|f - f_0\|_{L^2} > \omega_n(M_n \epsilon_n) \bigr\} .
\]
The posterior probability of the larger event tends to zero by
\eqref{eq:conclusion}.
\end{proof}
\begin{definition} (see \cite{knapik2011})
\label{def:credible_set}
    Denote the mean of the posterior $\widetilde{\Pi}$ by $v_{\rm map}$. Then the credible ball centred at $v_{\rm map}$ is defined as
    \begin{equation}
    \label{eq:credible_set}
        v_{\rm map} + B(r_{n,c}) := \{v \in V :||v - v_{\rm map}||_{V} < r_{n,c}  \}
    \end{equation}
    where $B(r_{n,c})$ is the ball centred at $v_{\rm map}$ with radius $r_{n,c}$. The constant, $c \in (0,1)$, denotes the desired credible level of $1-c$. The radius, $r_{n,c}$, is chosen such that
    \begin{equation}
        \widetilde{\Pi}_{n, \tau_n} (v_{\rm map} + B(r_{n,c}) \mid \widetilde{Y}_n) = 1-c.
    \end{equation}
The coverage of \cref{eq:credible_set} is then defined as
\begin{equation}
\label{eq:coverage}
    \widetilde{\Pi}_{n, \tau_{n}}(v_0 \in v_{\rm map} + B(r_{n,c}) \mid \widetilde{Y}_n)
\end{equation}
\end{definition}

\begin{proposition}[Credible Sets and Coverage]
\label{proposition:credible}
Let $\widetilde{C}_n(\widetilde{Y}_n)$ be a credible set for $v$ in
model~(L) and let
$C_n(Y_n) := \{ f \colon \mathbb{L}u_f \in
\widetilde{C}_n(Y_n - \tilde{g}) \}$.
Then the credible levels of $\widetilde{C}_n(\widetilde{Y}_n)$ and
$C_n(Y_n)$ are equal, and so are their coverage levels at $v_0$ and
$f_0$ respectively.
If moreover \eqref{eq:modulus} holds on $V_n$,
$\widetilde{C}_n(\widetilde{Y}_n) \subseteq V_n$, and
$\widetilde{C}_n(\widetilde{Y}_n) = \{ v \colon \|v - \hat{v}_n\|
\leq r_n \}$, then on the event $\{ v_0 \in
\widetilde{C}_n(\widetilde{Y}_n) \}$,
\[
  \operatorname{diam}_{L^2}\bigl( C_n(Y_n) \bigr)
  \;\leq\; 2\, \omega_n( 2 r_n ) .
\]
\end{proposition}
 
\begin{proof}
The first two assertions are \cite[proposition 2.2]{koers2024}. For
the third, the triangle inequality in $L^2$ and \eqref{eq:modulus}
give
\[
  \operatorname{diam}_{L^2}\bigl( C_n(Y_n) \bigr)
  = \sup_{v, w \in \widetilde{C}_n} \|e(v) - e(w)\|_{L^2}
  \leq 2 \sup_{v \in \widetilde{C}_n} \|e(v) - e(v_0)\|_{L^2}
  \leq 2\, \omega_n\Bigl( \sup_{v \in \widetilde{C}_n}
    \|v - v_0\| \Bigr),
\]
using that $\omega_n$ is non-decreasing. On the event
$\{ v_0 \in \widetilde{C}_n \}$ we have $\|\hat{v}_n - v_0\| \leq
r_n$, so $\|v - v_0\| \leq \|v - \hat{v}_n\| + \|\hat{v}_n - v_0\|
\leq 2 r_n$ for every $v \in \widetilde{C}_n$.
\end{proof}
 
\begin{remark}[Examples]
\label{rem:examples}
For $\omega_n(t) = L t$ one recovers
\cite[proposition 2.1]{koers2024} with rate $\epsilon_n$, as this is just the assumption that $e$ is Lipschitz. For a
H\"older modulus $\omega_n(t) = L t^{\kappa}$ with $\kappa \in (0,1]$
the rate is $\epsilon_n^{\kappa}$, since then $\rho_n = L M_n^{\kappa}
\epsilon_n^{\kappa}$ and $M_n^{\kappa} \to \infty$. 
\end{remark}
The following lemmas provide the crucial a-priori for our examples. 

\begin{lemma}\label{lem:c_continuous}
Let $\mathcal{K} \in L(V,U)$, and let $u_{f_{0}} = \mathcal{K}v_{f_{0}}$. 
Let $B_s(f_0)$ and $B_r(u_{f_0})$ be open balls around the true parameter $f_0$ and the corresponding solution respectively.
Let the operator $C(\cdot; f) \colon B_r(u_{f_0}) \to V$ satisfy a modulus of continuity $\omega_C$ uniformly for every $f \in B_s(f_0)$, 
and assume $C(\mathcal{K} v_{f_0}; f_0) = v_{f_0} - h$.
Assume there exists an operator $A \in L(\mathcal{F},V)$ that is invertible with bounded inverse and satisfies 
\begin{equation}\label{eq:H1}
  \| C(u; f_1) - C(u; f_2) - A(f_1-f_2)\|_{V}
  \leq \frac{\theta}{\|A^{-1}\|}\,\| f_1-f_2 \|_{\mathcal{F}}
\end{equation}
for every $u \in B_r(u_{f_0})$, $f_1,f_2 \in B_s(f_0)$ and a fixed constant $\theta \in (0,1)$.

Then there exist $\tilde{r}>0$ and a map 
$e\colon B_{\tilde{r}}(v_{f_0}) \to B_s(f_{0})$ with $e(v_{f_{0}})=f_{0}$ satisfying
\begin{align}
  \label{eq:HGsol}
  C(\mathcal{K} v_{f}; e(v_{f})) &= v_{f}-h,\\
\label{eq:HGmod}
  \| e(v_{1})-e(v_{2})\|_{\mathcal{F}}
  &\leq \tilde{\omega}(\| v_{1}-v_{2}\|_V),
\end{align}
where $\tilde{\omega}(t) := \frac{\|A^{-1}\|}{1-\theta}\left(t + \omega_C(\|\mathcal{K}\| t)\right)$,
for every $v_{f}, v_{1}, v_{2} \in B_{\tilde{r}}(v_{f_0})$.
\end{lemma}
\begin{proof}
  Since $\mathcal{K} \in L(V,U)$ is continuous, there exists $\hat{r}>0$ such that $\mathcal{K} B_{\hat{r}}(v_{f_0}) \subseteq B_{r}(u_{f_0})$.
  We define the mapping
  \begin{equation*}
    \Phi \colon \begin{cases}
       B_{\hat{r}}(v_{f_0}) \times B_s(f_0) \to \mathcal{F}\\ 
      (v,f) \mapsto A^{-1}C(\mathcal{K} v; f) - A^{-1}v + A^{-1}h.
    \end{cases}
  \end{equation*}
  We aim to apply the implicit function theorem (\cref{thm:graves}).
  We observe that by assumption $\Phi(v_{f_0},f_{0})=0$ at the reference point.
  We next prove that $\Phi(\cdot,f)$ satisfies a modulus of continuity i.e.
  \begin{align*}
  \|\Phi(v_1,f) - \Phi(v_2,f)\|_{\mathcal{F}} &\leq \|A^{-1}\| \| v_1 - v_2 + C(\mathcal{K} v_1; f) - C(\mathcal{K} v_2; f) \|_V \\
                                       &\leq \|A^{-1}\| \left(\|v_1-v_2\|_V + \omega_C(\|\mathcal{K} (v_1-v_2)\|_U)\right)\\
                                       &\leq \|A^{-1}\| \left(\|v_1-v_2\|_V + \omega_C(\|\mathcal{K}\| \|v_1-v_2\|_V)\right)\\
                                       &= (1-\theta) \tilde{\omega} (\|v_1-v_2\|_V).
  \end{align*}
  %Here $(1-\theta)\tilde{\omega}$ is a valid modulus of continuity as it is a non-decreasing function vanishing at zero. 
  %This modulus also does not depend on $f$, because $\omega_C$ is independent of $f$ by assumption.

 Subsequently, we compute the linear approximation with respect to $f$ by choosing $\tilde{A} = \mathrm{id}_\mathcal{F}$, to obtain 
  \begin{align*}
    \|\Phi(v,f_1) - \Phi(v,f_2) - \tilde{A}(f_1-f_2)\|_{\mathcal{F}} 
    &= \|A^{-1}C(\mathcal{K} v; f_1) - A^{-1}C(\mathcal{K} v; f_2)  - (f_1-f_2)\|_{\mathcal{F}}\\
    &\leq \|A^{-1}\| \|C(\mathcal{K} v; f_1) - C(\mathcal{K} v; f_2)  - A (f_1-f_2)\|_V\\
    &\leq \|A^{-1}\| \frac{\theta}{\|A^{-1}\|} \|f_1-f_2\|_{\mathcal{F}}\\
    &= \theta \|f_1-f_2\|_{\mathcal{F}}.
  \end{align*}
  Because $\|\tilde{A}\| = \|\tilde{A}^{-1}\| = 1$ and $\theta \in (0,1)$, we can apply \cref{thm:graves}
  to obtain the existence of a neighborhood $B_{\tilde{r}}(v_{f_0})$ and a continuous map $e\colon B_{\tilde{r}}(v_{f_0}) \to B_s(f_0)$ such that $\Phi(v,e(v))=0$, which rearranges exactly to $v + C(\mathcal{K} v; e(v)) = h$, for all $v \in B_{\tilde{r}}(v_{f_0})$.
  The modulus of continuity of $e$ follows directly from \cref{cor:graves} and by 
  \begin{align*}
    \|e(v_1)-e(v_2)\|_{\mathcal{F}} &\leq \frac{1}{1-\theta} (1-\theta) \tilde{\omega}(\|v_1-v_2\|_V) \\
                                  &= \frac{\|A^{-1}\|}{1-\theta}\left(\|v_1-v_2\|_V + \omega_C(\|\mathcal{K}\|\|v_1-v_2\|_V)\right) \\
                                  &= \tilde{\omega}(\|v_1-v_2\|_V).
  \end{align*}
\end{proof}

We want to apply \cref{thrm:tien_contraction_rate_early_stopping} in \cite{tienstra2025}, \cref{eq:linear_model}. For this analysis we need to work in the seuqence space setting. 
Particularly we will consider the observations, \cref{eq:linear_model} to be continuous observations.

\begin{claim}[Statistical equivalence of models for $d\ge1$, {\cite[Section 4]{reiss2008asymptotic}}]
Let $\mathbb{E}_{n,r}^{d}$ be the experiment obtained from observing independent random design points
$X_{i}\sim U_F([0,1]^{d})$, $i=1,\dots,n$, and the regression
\begin{equation}
  Y_{i}=f(X_{i})+\sigma\varepsilon_{i},\qquad i=1,\dots,n,
\end{equation}
for $n\in\mathbb{N}$ and $f:[0,1]^{d}\to\mathbb{R}$ in some class $\mathcal{F}^{d}\subseteq L^{2}([0,1]^{d})$,
with i.i.d.\ random variables $\varepsilon_{i}\sim\mathcal{N}(0,1)$ independent of the design.
Let $\mathbb{G}_{n}^{d}$ be the Gaussian white noise experiment given by observing
\begin{equation}\label{eq:cont_linear_observations}
  dY(x)=f(x)\,dx+\frac{\sigma}{\sqrt{n}}\,dW(x),
\end{equation}
where $dW$ is Gaussian white noise in $L^{2}L^{2}([0,1]^{d})$, understood weakly as the Gaussian process
$\psi\mapsto\langle\psi,dY\rangle$ indexed by $\psi\in L^{2}$ and realised on $L^2L^{2}([0,1]^{d})$.
Then, see \cite[Theorem 4.3]{reiss2008asymptotic}, for $d$-dimensional periodic Sobolev classes
$\mathcal{F}^{d}_{S,per}(s,R)$ with regularity $s>d/2$, the nonparametric regression experiment
$\mathbb{E}_{n,r}^{d}$ with random design and the Gaussian shift experiment $\mathbb{G}_{n}^{d}$ are
asymptotically equivalent as $n_{0},n\to\infty$ with $n_{0}=o(n^{1/2})$, and the Le Cam distance
satisfies
\begin{equation}\label{eq:lecam_distance_bound}
  \Delta_{\mathcal{F}^{d}_{S,per}(s,R)}\big(\mathbb{E}_{n,r}^{d},\mathbb{G}_{n}^{d}\big)
  \;\lesssim\; n^{-1/2}n_{0}+\sigma^{-1}R\,n_{0}^{1/2-s/d}.
\end{equation}
Moreover, on the event $(\psi_i)(X_j)_{i,j\leq n}$ is invertible, see \cite[Section 3.1]{reiss2008asymptotic}, and conditionally on the random design points, we can
equivalently observe
\begin{equation}\label{eq:sequence_space_obs}
  Z:=\sum_{i=1}^{n}\langle Y,\psi_{i}\rangle_{n}\,\psi_{i}
    =\Pi_{n}f+\frac{\sigma}{\sqrt{n}}\big(\Pi_{n}\!\mid_{S_{n}}\big)^{1/2}\xi \;\in\; S_{n},
\end{equation}
where $\xi$ is $\mathcal{N}(0,\mathrm{Id}_{S_{n}})$.
\end{claim}

\begin{proof}
The claim that \eqref{eq:lecam_distance_bound} holds is a direct application of
\cite[Theorem 4.3]{reiss2008asymptotic}; the optimal choice $n_{0}\sim n^{d/(2s+d)}$ yields
$\Delta\lesssim n^{(d-2s)/(2d+4s)}$ \cite[Remark 4.4]{reiss2008asymptotic}.
The claim that we can equivalently observe
\eqref{eq:sequence_space_obs} is \cite[Theorem 3.1]{reiss2008asymptotic}, applied conditionally on the
design as in \cite[Section 4.1]{reiss2008asymptotic}.
\end{proof}

To do this, we must check that the assumptions of \cref{thrm:tien_contraction_rate_early_stopping} are satisfied. For convenience, we collect all of the assumptions here. 
\begin{assumption}
We make the following assumptions :
    \label{assumption:1}
    \begin{enumerate}
        \item $\mathcal{K}\colon V \to U_F$ is a self-adjoint compact linear operator.
        \item There exist some $p >0$ such that for all $n \in \mathbb{N}$, the eigenvalues of $\mathcal{K}$ decay polynomially and
          \begin{equation}
         \kappa_i \asymp i^{-p/d} \quad i \ \forall n \in \mathbb{N}.
        \end{equation}
        \item The projection dimension is chosen as $D(n) \asymp n^{1/(2p + d)}$, where $P_n$ projects onto a $D(n)$-dimensional subspace of $V$ (see \cref{eq:projector}).
         \item The prior covariance operator $C_0$ is diagonalisable with respect to the basis from $\mathcal{K}^\ast\mathcal{K}$. And therefor $C_0$ and $\mathcal{K}$ commute.
         \item For a fixed ground truth $v_0 \in V$, with $v_0 \in H_d^{\beta'}$, define $ \widetilde{v}_0 := C_0^{-1/2}v_0$. We define $\beta$ to be such that $\widetilde{v}_0 \in H_d^\beta$, following \cref{eq:H_beta_d}, i.e. $\beta$ is the assumed regularity of the whitened truth $\widetilde{v}_0$ (which need not equal the true regularity $\beta'$ of $v_0$).
         \item The eigenvalues of $C_0$ with respect to the basis of $\mathcal{K}^*\mathcal{K}$ have the following structure (these are distinct from, and not to be confused with, the eigenvalues of $-\Delta$ used for the fixed, non-adaptive prior of \cref{eq:gaussian_prior})
        \begin{equation*}
             \lambda_i \asymp i^{-1-2\alpha/d} \quad i \in \mathbb{N}
         \end{equation*}
         for some $\alpha >0$ such that
         \begin{equation}
             \beta \leq d + 2\alpha + 2p
         \end{equation}
         holds.
         \item The stopping criterion \cref{eq:dp_stopping_time}, is chosen such that $\kappa \asymp D(n) /n$.
    \end{enumerate}
\end{assumption}
Suppose that the prior in \cref{eq:linearized_post}, $\widetilde{\Pi}_n(v_n)$ is now depending on a hyperparameter $\tau_n^2$. We will denote the prior now as  $\widetilde{\Pi}_{n, \tau_n}(v_n)$ to express the dependence on the hyperparameter $\tau_n$. The posterior given this prior is then also dependent on $\tau_n$ and will be denoted as
\begin{equation}
    \label{eq:linerised_post_tau}
    \widetilde{\Pi}_{n, \tau_n}( v \mid \widetilde{Y}_n).
\end{equation}
In the following lemma, we prove that the linearised posterior \cref{eq:linerised_post_tau}, contracts at the optimal rate for $\tilde{v}_0$ when $\tau_n = \tau_{\rm dp}$ \cref{eq:dp_stopping_time}.  
\begin{lemma}
\label{lem:linear_post_contracts}
   Let $v_0 \in H_d^{\beta^\prime}$. Let $\mathcal{K}$, be as in \cref{rem:spillting_pde}.  Fix the prior covariance such that \cref{assumption:1} hold. Suppose the prior is
   \begin{equation}
       \widetilde{\Pi}_{n, \tau_n}(v) \sim \mathcal{N}(0, \tau^2_{\rm dp}C_0)
    \end{equation}
  then 
   \begin{equation}
   \label{eq:v0_contraction_rate}
       \widetilde{\Pi}_{n,\tau_{\rm dp }} \left( \widehat{v}_{\tau_{\rm dp }} : || \widehat{v}_{\tau_{\rm dp }} - v_0||^2_{\ell^2(\mathbb{N})} \geq M_n \epsilon_n\right) \overset{ \mathbb{P}^{(n)}_{v_0}}{\rightarrow} 0.
    \end{equation}
for $\epsilon_n \asymp n^{-2\beta /( 2\beta + 2\alpha + 2p + 1 + d)}$, and $M_n \rightarrow \infty$.
\end{lemma}
\begin{proof}
    For $\mathcal{K}$, is as in \cref{rem:spillting_pde} assumption 1 in \cref{assumption:1} holds. As the rest of assumptions in \cref{assumption:1} holds, \cref{eq:v0_contraction_rate} follows directly from \cref{thrm:tien_contraction_rate_early_stopping}.
\end{proof}
We now show that this posterior contraction rate can be pulled back to \cref{eq:original_post}. To do this, we need to satisfy the conditions of Proposition 2.1 of \cite{koers2024}, see \cref{proposition:koers_contraction_pi_tilde} in the appendix for reference.

\begin{theorem}
\label{thrm:orginal_post_contracts}
 Suppose also that \cref{lem:linear_post_contracts} holds and that \cref{assumption:modulus} is satisfied. If $c$ is continuous in $u_f$, 
 then the posterior for $f$, 
    \begin{equation}
        \Pi_{n, \tau_{\rm dp}}(f_0 \mid Y_{n})
    \end{equation}
    contracts to $f_0$ at rate $\omega_n(\epsilon_n)$ on $V_n$, for $V_n$ as in \cref{assumption:modulus}. 

    If $c$ is differentiable in $u_f$, 
    then the posterior for $f$, 
    \begin{equation}
        \Pi_{n, \tau_{\rm dp}}(f_0 \mid Y_{n})
    \end{equation}
    contracts to $f_0$ at same rate $\epsilon_n$ on $V_n$, for $V_n$ as in \cref{assumption:modulus}. 
\end{theorem}
\begin{proof}
When $c$ is continuous in $u_f$ by \cref{lem:c_continuous} there exists a continuous map $e$ such that the conditions in \cref{proposition:transfer} hold. Similarly when $e$ is differentiable in $u_f$ by \cref{prop:nemytskii_modulus} there exists a Lipschitz map $e$ such that the conditions in \cref{proposition:koers_contraction_pi_tilde} hold. And thus the two claims following from applying the respective propositions. 
\end{proof}
As we are in the Bayesian setting, the posterior distribution provides a measure of uncertainty of our estimator. In the following results, we show that the posterior spread is a frequentist measure of uncertainty. 
\begin{lemma}
\label{lem:coverage_lin_post}
    For fixed $\alpha > 0$, if $v_{0,i} = C i^{-1-2\beta'/d}$ and $\widetilde{v}_0 = C i^{-1-2\beta/d}$ for all $i=1,...,D(n)$ and $\beta \leq d + 2 \alpha + 2p$, then as $n \rightarrow \infty$, $ \widetilde{\Pi}_{n, \tau_{\rm dp}}$  has frequentist coverage 1.
\end{lemma}
\begin{proof}
    This follows directly as a consequence of \cref{corollary:tien_assympotic_coverage} (see Lemma 2.8 in \cite{tienstra2025}). 
\end{proof}
We show that the coverage of \cref{eq:linearized_post} can be transferred to \cref{eq:original_post}. 
\begin{corollary}
\label{corollary:coverage_original_post}
    For fixed $\alpha > 0$, if $v_{0,i} = C i^{-1-2\beta'/d}$ and $\widetilde{v}_0 = C i^{-1-2\beta/d}$ for all $i=1,...,D(n)$ and $\beta \leq d + 2 \alpha + 2p$, then as $n \rightarrow \infty$, $ \Pi_{n, \tau_{\rm dp}}$  has frequentist coverage 1. Moreover if the credible set of $\widetilde{\Pi}(\cdot \mid \widetilde{Y}_n)$ haves radius $r_{n,c}$, then
    if $c$ is only continuous in $u_f$ then the radius of credible set in $\Pi_{n,\tau_{\rm dp}}(\cdot \mid Y_n)$ are of the order $\omega_n(r_n,c)$ by \cref{proposition:credible}. However, if $c$ is differentiable in $u_f$ then the radius of credible set in $\Pi_{n,\tau_{\rm dp}}(\cdot \mid Y_n)$ are of the same order.
\end{corollary}
\begin{proof}
    By \cref{proposition:koers_credible_sets_equal} the credible sets \cref{eq:credible_set} of $\widetilde{\Pi}_{n, \tau_{\rm dp}}(\cdot \mid \widetilde{Y})$ and $\Pi_{n, \tau_{\rm dp}}(\cdot \mid Y)$ centered at $v_{\rm map}$ and $f_{\rm map}$ respectively are the same. As $v_{\rm map} \rightarrow v_0$ and $f_{\rm map} \rightarrow f_0$ as $n \rightarrow \infty$, consequence of \cref{lem:coverage_lin_post} and \cref{thrm:orginal_post_contracts}, then asymptotically the coverage \cref{eq:coverage} of $\widetilde{\Pi}_{n, \tau_{\rm dp}}(\cdot \mid \widetilde{Y})$ and  $\Pi_{n, \tau_{\rm dp}}(\cdot \mid Y)$ is equal. Thus, the coverage of original posterior, \cref{eq:original_post}, $\Pi_{n,\tau_{\rm dp}}(\cdot \mid Y)$, which now depends on $\tau_{\rm dp}$ through  $\widetilde{\Pi}_{n, \tau_{\rm dp}}(\cdot \mid \widetilde{Y})$ via the solution map \cref{eq:solution_map_e}, is asymptotically 1, \cref{lem:coverage_lin_post}. Moreover, in the region where $e$, the solution map, is Lipschitz, the $r_{n,c}$ for both \cref{eq:linearized_post} and \eqref{eq:original_post} are of the same order.  Moreover in the case that $c$ is only continuous in $u_f$ then the radius of credible set in $\Pi_{n,\tau_{\rm dp}}(\cdot \mid Y_n)$ are of the order $\omega_n(r_n,c)$ by \cref{proposition:credible}.
\end{proof}

%%%%%%%%%%%%%%%%%%%%%%%%%%%%%%%%%%%%%%
\section{Extension to parabolic semi-linear PDEs} \label{sec:extensions}
%%%%%%%%%%%%%%%%%%%%%%%%%%%%%%%%%%%%%%
The theoretical framework presented in \cref{sec:main_results} assumed that the linearised inverse operator $\mathcal{K}$ is self-adjoint (see \cref{assumption:1}). 
This naturally covers stationary elliptic equations where the spatial differential operator $\mathbb{L}$ is self-adjoint. 
However, for time-dependent parabolic problems, the forward space-time differential operator is no longer self-adjoint. 
In this section, we outline how our early-stopping and linearisation methodology can be extended to space-time domains using the Singular Value Decomposition (SVD) of compact operators, and how the abstract implicit function theory applies seamlessly to space-time function spaces \cite{ladyzenskaja_linear_1968}.

\subsection{Linearisation and SVD in Space-Time}
Let $\mathcal{O}_T := (0, T) \times \mathcal{O}$ be a space-time cylinder. We generalize our framework to semilinear parabolic equations of the form:
\begin{equation}
\label{eq:parabolic_pde}
    \left\{ 
    \begin{alignedat}{4}
        \partial_t u_f + \mathbb{L} u_f - C(u_f; f) &= \, h \quad &&\text{in } \mathcal{O}_T, \\
        B_j u_f &= \, g_j \quad &&\text{on } \partial\mathcal{O} \times (0,T), \quad j = 0, \dots, m-1, \\
        u_f(0, x) &= \, u_0(x) \quad &&\text{on } \mathcal{O},
    \end{alignedat}
    \right.
\end{equation}
where $\mathbb{L}$ is a symmetric, strictly elliptic spatial differential operator of order $2m$ with strictly positive eigenvalues (e.g., for $2m=2$, $\mathbb{L} = -\Delta$). The operators $\{B_j\}_{j=0}^{m-1}$ form a normal system of boundary trace operators, e.g. $B_j = \partial_\nu^j$ for standard Dirichlet boundary conditions.
Defining the full parabolic operator $\mathcal{P} := \partial_t + \mathbb{L}$, we again linearise the problem by isolating the non-linear terms. 
Assuming the boundary data $g_j$ and initial data $u_0$ are sufficiently smooth and compatible, we can construct a known lifting function $\tilde{g}$ satisfying the initial and boundary conditions. 
We define our linearised state variable as $v = \mathcal{P} (u_f - \tilde{g})$. 
The inverse operator $\mathcal{K} = \mathcal{P}^{-1}$ solves the forward parabolic equation with homogeneous initial and boundary conditions, mapping a space-time source term $v$ to the solution displacement $u_f - \tilde{g}$. 
The adjoint operator $\mathcal{K}^* = (-\partial_t + \mathbb{L})^{-1}$ solves the \textit{backward} parabolic equation. 
Consequently, $\mathcal{K} \neq \mathcal{K}^*$, and they do not share an eigenbasis.
To recover the diagonal sequence-space model \cref{eq:truncated_model}, we rely on the singular system of $\mathcal{K}$. 
Because $\mathcal{K}$ remains a compact operator from $V := L^2(\mathcal{O}_T) \to L^2(\mathcal{O}_T)$, it admits a singular value decomposition (see e.g. \cite{kirsch_introduction_2021}) by $(\sigma_i, \phi_i, \psi_i)_{i \in \mathbb{N}}$, where
\begin{equation}
    \mathcal{K} \phi_i = \sigma_i \psi_i, \qquad \mathcal{K}^* \psi_i = \sigma_i \phi_i.
\end{equation}
Here, $(\phi_i)_{i \in \mathbb{N}}$ and $(\psi_i)_{i \in \mathbb{N}}$ each form an orthonormal basis for the domain and range of $\mathcal{K}$ respectively. 
By requiring the prior covariance $C_0$ to be diagonalisable with respect to the singular vectors $(\phi_i)$, and by projecting the observations $\widetilde{Y}_n = \mathcal{K}v + n^{-1/2} \Xi$ onto the singular vectors $(\psi_i)$, we compute the inner products:
\begin{align*}
    \langle \widetilde{Y}_n, \psi_i \rangle_{L^2} &= \langle \mathcal{K}v, \psi_i \rangle_{L^2} + n^{-1/2} \langle \Xi, \psi_i \rangle_{L^2} \\
    &= \langle v, \mathcal{K}^* \psi_i \rangle_{L^2} + n^{-1/2} \xi_i \\
    &= \sigma_i \langle v, \phi_i \rangle_{L^2} + n^{-1/2} \xi_i.
\end{align*}
By defining $\widetilde{Y}_{n, i} := \langle \widetilde{Y}_n, \psi_i \rangle$ and $v_i := \langle v, \phi_i \rangle$, we perfectly recover the diagonal sequence-space model
\begin{equation}
    \widetilde{Y}_{n, i} = \sigma_i v_i + n^{-1/2} \xi_i, \quad i = 1, \dots, D(n).
\end{equation}
Because the noise coordinates $\xi_i = \langle \Xi, \psi_i \rangle$ remain independent standard normal variables, the statistical structure of the linearised problem is identical to the self-adjoint case. 

\subsection{The Nemytskii Operator and Solution Map}
To pull the linearised contraction rates back to the original parameter $f$, we rely on the existence of a Lipschitz solution map $e(v) = f$. A major strength of our general theory in \cref{sec:main_results} is that \cref{ass:nonlinearity} and the Implicit Function Theorem application (\cref{lem:c_continuous}) are formulated entirely in abstract Banach spaces. 

For the time-dependent case, we define the base space as $V = L^2(\mathcal{O}_T)$. Because $\mathbb{L}$ has strictly positive eigenvalues, standard parabolic regularity guarantees that the operator $\mathcal{K} = (\partial_t + \mathbb{L})^{-1}$ is a bounded linear map from $V$ into the parabolic Sobolev space 
\begin{equation}
U = L^2\big(0,T; H^{2m}(\mathcal{O})\big) \cap H^1\big(0,T; L^2(\mathcal{O})\big).
\end{equation}

\begin{proposition}
\label{prop:parabolic_solution_map}
Let $\mathcal{O}_T = (0,T) \times \mathcal{O}$ be a space-time cylinder where $\mathcal{O} \subset \mathbb{R}^d$. 
Define the base space $V = L^2(\mathcal{O}_T)$ and the parabolic state space $U = L^2(0,T; H^{2m}(\mathcal{O})) \cap H^1(0,T; L^2(\mathcal{O}))$. 
Let $\mathcal{K} = (\partial_t + \mathbb{L})^{-1} \in L(V,U)$. 

Suppose the time-dependent non-linearity $c: \mathcal{O}_T \times \mathbb{R}^M \times \mathbb{R}^{M_f} \to \mathbb{R}$ satisfies the Carathéodory conditions (C1)-(C3) of \cref{ass:nonlinearity} almost everywhere on $\mathcal{O}_T$, subject to the following dimensional adjustments:
\begin{enumerate}[label=(\roman*)]
  \item The spatial dimension $d$ bounding the polynomial growth exponents $q_\alpha$ and $p_\alpha$ in (C1) is replaced by the effective parabolic dimension $d_{\mathrm{eff}} = d + 2m$.
    \item The local and global multi-index sets $A_{\mathrm{loc}}$ and $A_{\mathrm{glob}}$ are defined with respect to the continuous embeddings of the parabolic Sobolev space $U$ over $\mathcal{O}_T$.
\end{enumerate}
Then the induced space-time Nemytskii operator $C: U \times \mathcal{F} \to V$ is well-defined and satisfies the abstract mapping properties of \cref{prop:nemytskii_modulus}. 
Consequently, there exists a neighborhood $B_{\tilde{r}}(v_{f_0}) \subset V$ and a solution map $e: B_{\tilde{r}}(v_{f_0}) \to \mathcal{F}$ such that $e(\mathcal{P}(u_f - \tilde{g})) = f$, which admits a modulus of continuity.
\end{proposition}
\begin{proof}
  Because the differential operator $\partial_t + \mathbb{L}$ is strictly parabolic with a spatial operator of order $2m$, the anisotropic Sobolev embedding on the cylinder $\mathcal{O}_T$ scales such that time acts as $2m$ spatial dimensions \cite{ladyzenskaja_linear_1968}.
Thus, the effective embedding dimension is $d_{\mathrm{eff}} = d + 2m$. 
Under this adjusted dimension, the continuous embeddings $U \hookrightarrow L^q(\mathcal{O}_T)$ dictate the sets $A_{\mathrm{loc}}$ (jets that embed into $L^\infty(\mathcal{O}_T)$) and $A_{\mathrm{glob}}$ (jets that embed only into Lebesgue spaces). 
By assumption, the time-dependent non-linearity $c(t,x,\xi;\zeta)$ respects the growth exponents constrained by $d_{\mathrm{eff}}$. Consequently, the identical Hölder and Jensen inequality arguments used in the proof of \cref{prop:nemytskii_modulus} apply on $L^2(\mathcal{O}_T)$ instead of $L^2(\mathcal{O})$. 
This guarantees the space-time Nemytskii operator maps $U \times \mathcal{F} \to V$ and admits a modulus of continuity $\omega_C$. Because $\mathcal{K} \in L(V,U)$ is a bounded linear map by standard parabolic regularity, and condition (C3) ensures the Fréchet derivative with respect to $f$ has a bounded inverse, all conditions of the abstract Implicit Function Theorem (\cref{lem:c_continuous}) are satisfied on the Banach spaces $V$ and $U$. The existence and local Lipschitz continuity (or modulus of continuity) of the inverse map $e: v \mapsto f$ immediately follows.
\end{proof}

%Because \cref{lem:c_continuous} depends only on the mapping properties between $V$ and $U$, the entire deterministic pullback machinery holds unchanged. 
%One only needs to verify that the time-dependent non-linearity $c(t, x, J^k u; J^l f)$ satisfies the Carathéodory conditions of \cref{ass:nonlinearity} almost everywhere on the cylinder $\mathcal{O}_T$, with the growth bound in (C1) applied using the effective space-time dimension $d_{\mathrm{eff}} = d+2$ in place of the spatial dimension $d$, reflecting the parabolic Sobolev embedding on $\mathcal{O}_T$ rather than the elliptic embedding on $\mathcal{O}$. 
%If these abstract polynomial bounds hold, the implicit function theorem immediately guarantees a locally Lipschitz solution map $e: V \to \mathcal{F}$. 

Utilizing \cref{prop:parabolic_solution_map}, the theoretical guarantees and asymptotic frequentist coverage can be transferred seamlessly to time-dependent inverse problems.

%%%%%%%%%%%%%%%%%%%%%%%%%%%%%%%%%%%%%%

\section{Examples: Theoretical Results} \label{sec:examples}

%%%%%%%%%%%%%%%%%%%%%%%%%%%%%%%%%%%%%%

In this section we demonstrate how the general results apply to a variety of examples. We begin with the canonical nonlinear example: the time-independent Schrödinger equation. In this example we also provide numerical results, see \cref{sec:enkbf}.

\subsection{Theoretical Results for Schrödinger}\label{subsec:results_schro}

In this section, we consider the canonical nonlinear example: the time-independent Schrödinger equation. This non-linear inverse problem is widely studied in the Bayesian inverse problems literature; see \cite{NicklWang, NicklLectureNotes2023, koers2024, giordano2020consistency}, among others. We demonstrate in detail that the results established in \cref{sec:main_results} apply. Additional examples where our theoretical framework is applicable can be found in \cite{koers2024}.

We are interested in the Bayesian problem of estimating  $f \in L^2(\mathcal{O})$ that is strictly positive $f_0 >0$ from noisy observations \cref{eq:linear_model}. To do this, we will apply the results of \cref{sec:main_results}. We give the problem in detail below, which can be originally found in the sources already mentioned. 

Let $\mathcal{O} \subseteq \mathbb{R}^d$, for $d \leq 2$. Let $f_0$ denote the ground truth. Suppose there exists a solution map $f=e(\mathbb{L}u_{f})$ \cref{eq:solution_map_e}, that is Lipschitz around $f_0$. Assume also that $u_f=u_{f_0}$ is the solution to
    \begin{equation}
    \label{eq:schrodinger}
        \begin{cases}
       -\Delta u_f +  f_0 u_f =0 \quad \text{on } \mathcal{O} \\
        u_f = g \quad \text{on } \partial \mathcal{O}. 
    \end{cases}
    \end{equation}
    where $g: \partial \mathcal{O} \rightarrow \mathbb{R}$ and is fixed. Writing this problem as a regression problem following \cref{sec:intro}, we also assume that $u_{f_0}$ is such that 
    \begin{align}
    \label{eq:regression_schro}
        Y_i &= G(f_0)(X_i) + \xi_i \\
            &= u_{f_0}(X_i) + \xi_i.
    \end{align}
In this problem $\mathbb{L}$ from \cref{eq:generic_pde} is $-\Delta$, where we take the negative Laplacian so that the eigenvalues are positive and the sign is preserved.
The log-likelihood over $f$ given the observations is
    \begin{equation}
    \label{likelihood}
        \ell_N(\theta) = -\frac{1}{2} \sum^{n}_{i=1} [ Y_i -\mathcal{G}(f)(X_i) ]^2, \quad f \in \mathbb{R}^{D(n)}.
    \end{equation}
    Following \cite{NicklWang}, we construct a Gaussian prior from the eigenvalues of the Laplacian. Denote the eigenvalues of the Laplacian, by $(\lambda_k)_{k \in \mathbb{N}}$. We must fix $\alpha$, and let $\tau_n$ be a scalar.  The prior for over $f$ is then
        \begin{equation}
             \label{eq:prior_nonlinear}
             \mathcal{N}(0,  \tau_n^2 \Lambda_\alpha ^{-1})
        \end{equation}
        where $\Lambda_\alpha = diag(\lambda^\alpha_1, ..., \lambda^\alpha_D)$.
        The posterior given observations $D^{(n)}= ((Y_1, X_1),...,(Y_n,X_n))$ is then
    \begin{align}
    \label{eq:post_schro}
        \Pi_{\tau_{\rm dp},n}(f \mid D^{(n)}) & \ \propto \ e^{\ell_N(\theta)} \, \Pi_n(f) \\
        & \ \propto \text{ exp}\left\{- \sum_{i=1}^n (Y_i - \mathcal{G}(\theta)(X_i))^2 - \tau_n^2 ||f||_{\ell^2(\mathbb{N})}^2\right\}.
    \end{align}
    The point estimator for $f_0$ is then
    \begin{equation}
        \widehat{f}_{MAP} \in \underset{f \in \mathbb{R}^{D(n)}}{\text{arg max}} \Pi_n(f \mid D^{(n)}).
        \end{equation}

    \begin{lemma}[Existence of Lipschitz solution map]
        \label{lem:schro_lip_solution_map}
        To obtain a Lipschitz bound on $e$, we verify the assumptions of \cref{prop:nemytskii_modulus}.
        We consider the candidate for the Fréchet derivative $A_{u_{f},{f}}(\eta,\zeta)  = M_f \eta +  M_{u_f-\tilde{g}}\zeta$.
        Here $M_f\colon L^2(\mathcal{O})\rightarrow L^2(\mathcal{O})$ is the multiplication operator multiplying by $f\in \mathcal{F}$. It is bounded as $\mathcal{F}\subseteq L^\infty(\mathcal{O})$ and 
        $
            \| u_f f\| \leq \|f\|_{L^\infty} \|u_f\|.
        $
        For every $(u_f,f)\in L^2(\mathcal{O})\times \mathcal{F}$, $\eta \in L^2(\mathcal{O})$, $\zeta \in \mathcal{F}$, and $A_{u_f,f}$ is a linear, bounded map $L^2(\mathcal{O}) \times \mathcal{F} \rightarrow L^2(\mathcal{O})$ satisfying
        $$
        %\begin{aligned}
            \|c(u_f+\eta, f + \zeta)  - c(u_f, f) - A_{u_f,f}(\eta,\zeta)\|_{L^2(\mathcal{O})}
            %&= \| M_{f+\zeta}(u_f+\eta-\tilde{g}) - M_{f}(u_f-\tilde{g}) - M_f \eta -  M_{u_f-\tilde{g}}\zeta \|\\
            = \| \zeta \eta \|_{L^2(\mathcal{O})}
            \leq \|\zeta\|_{\mathcal{F}} \|\eta\|_{L^2(\mathcal{O})} \in o(\|(\zeta,\eta)\|).
        % \end{aligned}
        $$
        As $f\mapsto M_f$ and $u_f\mapsto M_{u_f-\tilde{g}}$ are continuous $(u_f,f) \mapsto A_{u_f,f}$ is continuous as composition of continuous functions and therefore $c$ is continuously Fréchet differentiable.
        Furthermore $\partial_f c_{u_f,f} = M_{u_f-\tilde{g}}$ (the partial Fréchet derivative) has bounded inverse for every $u_f \in {B_\varepsilon(\tilde{g})}^c \subset L^2(\mathcal{O})$.
        Finally  $\mathcal{K}\colon L^2(\mathcal{O})\rightarrow L^2(\mathcal{O})$ is a linear, self-adjoint compact operator \cite{evans_partial_2010}.
    \end{lemma}
    \begin{remark}
    \label{rem:sequence_space_project_obs}
        To apply the results of \cref{sec:main_results}, we will consider the continuous observations 
        \begin{equation}
            Y_n = u_{f_0} + \frac{1}{\sqrt{n}}\Xi
        \end{equation}
        where the size of the grid has gone to $\infty$, i.e.\ the number of design points $X_i$ in \cref{eq:linear_model}. We then transform these observations using the method described in \cref{sec:intro}(the statistical equivalence in the Claim above \cref{eq:cont_linear_observations}), and project the linearised observations into sequence space. We can then make sense of  \cref{assumption:1} theoretically. To apply \cref{lem:linear_post_contracts}, we project the observations in sequence space by $P_n$. The result is that we have an estimator $\widehat{v}_{\rm map} \in \mathbb{R}^{D(n)}$. We can then take this estimator and transform it back into function space and then apply the solution map $e$ to $\mathbb{L}\widehat{v}_{\rm map}^{\rm func}$ to have an estimate for $f_0$. Practically, in this step, we need to evaluate on a grid, i.e.\ the numerical discretisation of $\mathcal{O}$ used in \cref{sec:enkbf}, see \cref{fig:schrodinger}.
    \end{remark}
  
    \begin{lemma}
    \label{lem:inverse_of_laplace}
        There exists a $\mathcal{K}: L_2(\mathcal{O}) \rightarrow L_2(\mathcal{O})$, such that for \cref{eq:schrodinger} where $\mathbb{L} = -\Delta$, we have that $\mathcal{K} = -\Delta^{-1}$ is the solution to \cref{eq:generic_pde}, and is self-adjoint and linear operator. 
    \end{lemma}
    \begin{proof}
        \cite[Chapter 6]{evans_partial_2010}
    \end{proof}
   
    \begin{lemma}
    \label{lem:lip_solution_map}
        For the problem of inferring $f_0$ from solutions of \cref{eq:schroedinger_eq}, there exists a solution map $e$ that is Lipschitz on the set $B^V(v_{f_0}) \subseteq V$. Moreover there exists nested subsets $B^V_{n}(v_{f_0})$ of $B^V(v_{f_0})$ such that
        \begin{equation}
            \widetilde{\Pi}(B^V_{n}(v_{f_0})\mid \widetilde{Y}_n) \rightarrow 1, \quad \text{as } n \rightarrow \infty.
        \end{equation}
    \end{lemma}
    \begin{proof}
       From \cref{eq:schroedinger_eq} we know that \cref{prop:nemytskii_modulus} holds. Thus, we know that there exists an $e$ that is Lipschitz on the set $B^V(v_{f_0}) \subseteq V$, say $B^V(v_{f_0}) = B_{\tilde r}(v_{f_0})$ for some $\tilde r>0$. We now construct the required nested sets directly from \cref{lem:linear_post_contracts}. By \cref{eq:v0_contraction_rate}, for the rate $\epsilon_n \asymp n^{-2\beta/(2\beta+2\alpha+2p+1+d)}$ and every sequence $M_n \to \infty$,
       \begin{equation}
           \widetilde{\Pi}_{n,\tau_{\rm dp}}\left(v : \|v - v_{f_0}\|_V^2 \geq M_n \epsilon_n \mid \widetilde{Y}_n\right) \overset{P_{v_0}}{\longrightarrow} 0.
       \end{equation}
       Fix such a sequence $M_n \to \infty$ with $M_n \epsilon_n \to 0$, chosen (without loss of generality, replacing $M_n\epsilon_n$ by $\sup_{m \geq n} M_m \epsilon_m$ if necessary) so that $M_n \epsilon_n$ is non-increasing in $n$, and define
       \begin{equation}
       \label{eq:ball_of_lip_solution_map}
           B^V_n(v_{f_0}) := \left\{ v \in V : \|v - v_{f_0}\|_V^2 \leq M_n \epsilon_n \right\}.
       \end{equation}
       Since $M_n \epsilon_n \to 0$, there exists $N$ such that $M_n \epsilon_n < \tilde r^2$ for all $n \geq N$, so $B^V_n(v_{f_0}) \subseteq B^V(v_{f_0})$ for $n \geq N$; discarding the finitely many $n < N$, $(B_n^V(v_{f_0}))_n$ is a nested (decreasing, since $M_n\epsilon_n$ is non-increasing) sequence of subsets of $B^V(v_{f_0})$, on each of which $e$ is Lipschitz. Finally,
       \begin{equation}
           \widetilde{\Pi}\left(B^V_n(v_{f_0}) \mid \widetilde{Y}_n\right) = 1 - \widetilde{\Pi}_{n,\tau_{\rm dp}}\left(v : \|v - v_{f_0}\|_V^2 > M_n \epsilon_n \mid \widetilde{Y}_n\right) \rightarrow 1
       \end{equation}
       in $P_{v_0}$-probability, as required.
    \end{proof}

    \begin{corollary}
        \label{corollary:schrodinger_contraction_rate}
        Suppose that $v_0\in H_d^{\beta'}$ and define $\widetilde{v}_0 := C_0^{-1/2} v_0$. Assume that $\tilde{v}_0\in H_d^{\beta}$ and $C_0$ is such that the eigenvalues with respect to the basis from $\mathcal{K}$ decay as $i^{-1-2\alpha/d}$, where $\alpha$ is such that $\beta < d + 2\alpha +2p$ holds. Let $D(n)$ as in \cref{assumption:1}. Let $\kappa \asymp D(n)/n$. Assume $\tau_{\rm dp}$ is estimated from \cref{eq:dp_stopping_time}. Then the posterior
        \begin{equation}
            \Pi_{\tau_{\rm dp},n}(f \mid Y_n)
        \end{equation}
        of the Bayesian inverse problem \cref{eq:regression_schro} of estimating the potential $f$ arising from \cref{eq:schrodinger} contracts rate $\epsilon_n \asymp n^{-2\beta/(2\beta + 2\alpha + 2p + 1 + d)}$ to $f_0$.
    \end{corollary}
    \begin{proof}
        We want to apply \cref{lem:coverage_lin_post}. First from \cref{lem:inverse_of_laplace}, we can find a $\mathcal{K}: L_2(\mathcal{O}) \rightarrow L_2(\mathcal{O})$, such that $\mathcal{K}$ is the inverse of $\mathbb{L}=-\Delta$ and is compact self-adjoint, and linear, and \cref{eq:generic_pde} holds. We can then form observations 
        \begin{align}
            \label{eq:schro_obs}
            \widetilde{Y}_n = -\Delta^{-1}v + n^{-1/2}\Xi
        \end{align}
        where $\tilde{g}$ is such that \cref{eq:generic_pde_boundary} holds. Such a $\tilde{g}$ exists when $\mathcal{O}$ is regular, by standard elliptic PDE theory \cite[Chapter 6]{evans_partial_2010}. By \cref{lem:lip_solution_map} when $u_{f_0} = \mathcal{K}v_0 + \tilde{g}$ is such that $\underset{x \in \mathcal{O}}{\text{inf}} u_{f_0} \geq r_0 > 0$, then the solution map on a ball around $\mathcal{K}v_0$ of radius $r_0$ is Lipschtz for all points in the ball \cref{eq:ball_of_lip_solution_map} (the lower bound keeps $f = v/(2(\mathcal{K}v+\tilde{g}))$, \cref{eq:solution_map_e}, away from its singularity at $\mathcal{K}v + \tilde{g} = 0$). We also have that the eigenvalues of $\mathcal{K} = (-\Delta)^{-1}$ (equivalently, the inverse eigenvalues of $-\Delta$) on $\mathcal{O} \subset \mathbb{R}^d$ with homogeneous boundary conditions, are such that
        \begin{equation}
        \label{eq:decay_laplace_eigenvalues}
            \lambda_i \asymp i^{-2/d}.
        \end{equation}
        So $\lambda_i$  decay polynomial with power $2/d$. By construction, the prior and prior covariance operator satisfy conditions in \cref{assumption:1}. Thus, the conclusion of the corollary follows from applying \cref{lem:coverage_lin_post} to the linearised posterior \eqref{eq:linearized_post}, and then applying \cref{thrm:orginal_post_contracts}.  
    \end{proof}
    We also have a result for the asymptotic coverage of the data-dependent posterior. Specifically that 
    \begin{corollary}
        \label{corollary:schrodinger_coverage}
    Suppose (C1) and (C2) of \cref{assumption:1} hold true, and that $\beta' > \beta$. For fixed $\alpha > 0$ such that $\beta \leq d + 2\alpha + 2p$ holds, if $v_{0,i} = C_i i^{-1-2\beta'/d}$ and $\widetilde{v}_{0,i} = C_i i^{-1-2\beta/d}$ for all $i=1,...,D(n)$, then as $n \rightarrow \infty$, $ \Pi_{n, \tau_{\rm dp}}$  has frequentist coverage 1.
    \end{corollary}
   
    \begin{proof}
         Recall that we denoted the posterior arising from Gaussian prior $\widetilde{\Pi}_n(v_0)$, with likelihood $L(\widetilde{Y}_{i} \mid v)$ as $\widetilde{\Pi}_{n,\tau_{\rm dp}}(v_0 \mid \widetilde{Y}_{i})$. As $\mathcal{K}$ satisfies \cref{assumption:1}, so by \cref{lem:coverage_lin_post}
    \begin{equation*}
        \widetilde{\Pi}_{n,\tau_{\rm dp}}\left(v_0 \in v_{\rm map}+ B(\tilde{r}_{n,c})\right) \rightarrow 1
    \end{equation*}
    as $n \rightarrow \infty$. Further then by \cref{corollary:coverage_original_post}
     \begin{equation*}
        {\Pi}_{n,\tau_{\rm dp}}\left(f_0 \in f_{\rm map}+ B(r_{n,c})\right) \rightarrow 1
    \end{equation*}
    as $n \rightarrow \infty$. More over if the credible set $v_{\rm map}+ B(\tilde{r}_{n,c}) \subseteq V_n$, where $V_n$ is the set such that the solution map is Lipschitz \cref{eq:ball_of_lip_solution_map}, $\tilde{r}_{n,c}$ is the diameter of the ball in the norm $\|\cdot \|_V$ and $r_{n,c}$ and is the diameter of the ball in the norm $\|\cdot \|_{L^2}$ then
    \begin{equation*}
        \tilde{r}_{n,c} \asymp r_{n,c} 
    \end{equation*}
    in probability under $v_0$ as consequence of \cref{proposition:koers_credible_sets_equal}.
    \end{proof}

\subsection{Allen-Cahn Equation}\label{subsec:results_allencahn}
To demonstrate that our framework easily handles non-linearities with polynomial growth, we consider the parameter estimation of the source term in the stationary Allen-Cahn equation (a canonical reaction-diffusion model). Let $\mathcal{O} \subset \mathbb{R}^d$ ($d \le 3$) and consider:
\begin{equation}
\label{eq:allencahn}
    \left\{ 
    \begin{alignedat}{4}
        - \Delta u_f + u_f^3 - u_f &= \, f \quad &&\text{in } \mathcal{O}, \\
        u_f &= \, 0 \quad &&\text{on } \partial\mathcal{O}.
    \end{alignedat}
    \right.
\end{equation}
Here, the unknown parameter to infer is the source field $f \in L^2(\mathcal{O})$. 
This casts perfectly into the general framework of \cref{eq:generic_pde}. We identify the linear operator $\mathbb{L} = -\Delta$ with inverse $\mathcal{K} = (-\Delta)^{-1}$, and the Nemytskii operator is given pointwise by $c(x, u; f) = -u^3 + u + f$. 
\begin{lemma}[Lipschitz solution map for Allen-Cahn]
    The non-linearity $c(x, u; f) = -u^3 + u + f$ satisfies conditions (C1)-(C3) of \cref{ass:nonlinearity}. Consequently, the inverse map $e(v) = f$ exists globally and is Lipschitz on bounded sets in $L^2(\mathcal{O})$.
\end{lemma}
\begin{proof}
    Because $d \le 3$, the space $U = H^2(\mathcal{O}) \cap H^1_0(\mathcal{O})$ embeds continuously into $L^\infty(\mathcal{O})$. 
    Therefore, the non-linearity is entirely local i.e. $A_{\mathrm{glob}} = \emptyset$. 
    Condition (C1) holds since $|-\xi^3 + \xi + \zeta| \le |\zeta| + |\xi| + |\xi|^3$, which is bounded by a constant for bounded inputs $|\xi| \le R$. 
    Condition (C2) holds since the polynomial is locally Lipschitz
    \begin{equation*}
        |(-\xi^3 + \xi + \zeta) - (-\eta^3 + \eta + \zeta)| \le |\xi - \eta| \left(1 + |\xi^2 + \xi\eta + \eta^2|\right) \le L_R |\xi - \eta|.
    \end{equation*}
    For (C3), the partial Fréchet derivative with respect to $f$ is simply the identity operator $A = \mathrm{id}$, which obviously has a bounded inverse. 
    
    Applying the linearisation method, we set $v = -\Delta u_f$, so $u_f = \mathcal{K}v$. The PDE rearranges directly to an explicit solution map:
    \begin{equation}
        f = e(v) = v + (\mathcal{K}v)^3 - \mathcal{K}v.
    \end{equation}
    For any $v_1, v_2$ in a bounded ball $V_n \subset L^2(\mathcal{O})$, the Sobolev embedding ensures $\|\mathcal{K}v_i\|_{L^\infty} \le C$. Evaluating the difference gives:
    \begin{align*}
        \|e(v_1) - e(v_2)\|_{L^2} &\le \|v_1 - v_2\|_{L^2} + \|\mathcal{K}(v_1 - v_2)\|_{L^2} + \|(\mathcal{K}v_1)^3 - (\mathcal{K}v_2)^3\|_{L^2} \\
        &\le \left( 1 + \|\mathcal{K}\| + C' \|\mathcal{K}\| \right) \|v_1 - v_2\|_{L^2}.
    \end{align*}
    Thus, $e(v)$ is Lipschitz on $V_n$.
\end{proof}

\begin{corollary}[Contraction rate and Coverage for Allen-Cahn]
\label{cor:allencahn_contraction}
Let the true source $f_0$ induce a true linearised parameter $v_0 \in H_d^{\beta'}$. Under the identical prior and stopping-time assumptions as \cref{corollary:schrodinger_contraction_rate}, the early-stopped posterior $\Pi_{\tau_{\rm dp},n}(f \mid Y_n)$ contracts to $f_0$ at the optimal rate $\epsilon_n$ for $\tilde{v}_0$ and yields asymptotic frequentist coverage of $1$.
\end{corollary}
\begin{proof}
We only prove that there exsists a Lipschitz map between $v$ and $f$. Because $e(v)$ is explicit and Lipschitz on bounded sets, it admits a modulus of continuity $\omega_n(t) \propto t$. The rest of the proof follows identically as in \cref{corollary:schrodinger_contraction_rate} and \cref{corollary:schrodinger_coverage} by applying \cref{thrm:orginal_post_contracts} and \cref{corollary:coverage_original_post}, the optimal contraction and valid frequentist coverage established for the linearised sequence space immediately pull back to the posterior over the source term $f_0$.
\end{proof}

\subsection{Parabolic Allen-Cahn Equation}
To make the abstract time-dependent framework concrete, we consider the inverse source problem for the time-dependent Allen-Cahn equation. 
Let $\mathcal{O} \subset \mathbb{R}^3$ be a bounded domain with a $C^2$ boundary. 
We aim to recover a space-time varying source parameter $f \in L^2(\mathcal{O}_T)$ from noisy observations of the linearised state. 
The governing equation on the space-time cylinder $\mathcal{O}_T := (0,T) \times \mathcal{O}$ is:
\begin{equation}
\label{eq:parabolic_ac}
    \left\{ 
    \begin{alignedat}{4}
        \partial_t u_f - \Delta u_f + u_f^3 - u_f &= \, f \quad &&\text{in } \mathcal{O}_T, \\
        u_f &= \, 0 \quad &&\text{on } \partial\mathcal{O} \times (0,T), \\
        u_f(0, x) &= \, 0 \quad &&\text{on } \mathcal{O}.
    \end{alignedat}
    \right.
\end{equation}
Here, the spatial operator is $\mathbb{L} = -\Delta$ (order $2m = 2$), and the non-linearity is given pointwise by $c(t, x, \xi; \zeta) = -\xi^3 + \xi + \zeta$. The highest jet order for both the state and the parameter is $k = 0$ and $l = 0$. 
We define the base space $V = L^2(\mathcal{O}_T)$ and the parameter space $\mathcal{F} = L^2(\mathcal{O}_T)$. By standard parabolic regularity, the linearised forward operator $\mathcal{K} = (\partial_t - \Delta)^{-1}$ is a bounded isomorphism from $V$ onto the parabolic Sobolev space:
\begin{equation}
    U = L^2\big(0,T; H^2(\mathcal{O}) \cap H^1_0(\mathcal{O})\big) \cap H^1\big(0,T; L^2(\mathcal{O})\big).
\end{equation}

\begin{proposition}
\label{prop:allen_cahn_spacetime}
The Nemytskii operator $C(u;f) = -u^3 + u + f$ satisfies the abstract Carathéodory growth and continuity conditions (C1)-(C3) as mappings from $U \times \mathcal{F} \to V$. Consequently, the solution map $e(v) = f$ is locally Lipschitz.
\end{proposition}
\begin{proof}
For space-time cylinders, the effective scaling dimension is $d_{\text{eff}} = d + 2$. For $d=3$, $d_{\text{eff}} = 5$. 
The parabolic Sobolev embedding theorem (see, e.g., \cite{ladyzenskaja_linear_1968}) yields the continuous continuous embedding $U \hookrightarrow L^{\frac{2 d_{\text{eff}}}{d_{\text{eff}} - 4}}(\mathcal{O}_T) = L^{10}(\mathcal{O}_T)$. 
Because $U$ does not embed into $L^\infty(\mathcal{O}_T)$, the state variable index $\alpha = 0$ strictly belongs to the global index set $A_{\mathrm{glob}}$ and must be bounded using Lebesgue norms. 

The growth of the Nemytskii operator is dominated by the cubic term. From the triangle inequality, $|c(t, x, \xi; \zeta)| \le |\xi|^3 + |\xi| + |\zeta|$. 
To satisfy (C1), we require exponents $q_0 = 3$ for the state and $p_0 = 1$ for the parameter, ensuring the operator maps into $L^2(\mathcal{O}_T)$. 
This requires the state space to embed into $L^{2q_0}(\mathcal{O}_T) = L^6(\mathcal{O}_T)$. 
Because the parabolic embedding provides $U \hookrightarrow L^{10}(\mathcal{O}_T)$ and $10 \ge 6$, the required norm is continuously bounded i.e. $\|u^3\|_{L^2} = \|u\|_{L^6}^3 \lesssim \|u\|_U^3$. 
Since $\mathcal{F} = L^2(\mathcal{O}_T)$, the parameter bound for $p_0=1$ holds trivially. 
Thus, (C1) is satisfied.

We factorise the difference of the non-linearity:
\begin{equation}
    |c(t,x,\xi_1;\zeta_1) - c(t,x,\xi_2;\zeta_2)| \leq \frac{3}{2} \left( |\xi_1|^2 + |\xi_2|^2 \right) |\xi_1 - \xi_2| + |\xi_1 - \xi_2| + |\zeta_1 - \zeta_2|.
\end{equation}
This perfectly matches the abstract structural requirement of \cref{eq:continuity_bound} in (C2) with $\gamma_0 = 1$ and $q_0 = 3$. 
Evaluating the $L^2(\mathcal{O}_T)$ norm of the cubic difference term allows us to apply Hölder's inequality with conjugate exponents $p = 3$ and $q = \frac{3}{2}$
\begin{equation}
    \left\| (\xi_1^2 + \xi_2^2)(\xi_1 - \xi_2) \right\|_{L^2(\mathcal{O}_T)} 
    \leq 
    \left\| \xi_1^2 + \xi_2^2 \right\|_{L^3(\mathcal{O}_T)} \left\| \xi_1 - \xi_2 \right\|_{L^6(\mathcal{O}_T)}.
\end{equation}
Because $\| \xi^2 \|_{L^3} = \| \xi \|_{L^6}^2$, all terms are bounded within the $L^6(\mathcal{O}_T)$ embedding of $U$. 
This proves the operator is locally Lipschitz in $U$, satisfying (C2).

The Fréchet derivative of $C(u; f)$ with respect to the parameter $f$ in the direction of a perturbation $\delta f$ is strictly linear $A(\delta f) = \delta f$. 
Therefore, the operator $A$ is the identity map. 
Its inverse $A^{-1}$ trivially maps $V = L^2(\mathcal{O}_T)$ to $\mathcal{F} = L^2(\mathcal{O}_T)$ with bounded norm $\|A^{-1}\|_{\mathcal{L}(V,\mathcal{F})} = 1$. 
The exact Taylor remainder of the linearisation with respect to $f$ is strictly zero, meaning \cref{eq:param_approx} holds with $\theta = 0$. 

By the abstract Implicit Function Theorem formulation (\cref{lem:c_continuous}), the non-linearity assumptions (C1)-(C3) guarantee the solution map $e: v \mapsto f$ is locally Lipschitz. 
\end{proof}

\subsection{Darcy Flow}\label{subsec:results_darcy}

We consider the standard formulation of Darcy flow, which models the pressure field $u$ of a fluid moving through a porous medium with an unknown permeability field $f$. Let $\mathcal{O} \subset \mathbb{R}^d$ ($d \le 3$) be a bounded domain with a $C^2$ boundary $\partial \mathcal{O}$ and outward unit normal $\nu$. The governing equation is
\begin{equation}\label{eq:darcy_standard}
    \left\{ 
    \begin{alignedat}{4}
        \nabla \cdot (f \nabla u) &= \, h \quad &&\text{in } \mathcal{O}, \\
        u &= \, g \quad &&\text{on } \partial\mathcal{O}.
    \end{alignedat}
    \right.
\end{equation}
Unlike the semilinear examples treated in \cref{sec:main_results}, the highest-order derivative of $u$ in \cref{eq:darcy_standard} is multiplied by the unknown parameter $f$, making the equation quasilinear. Consequently, the abstract Implicit Function Theorem of \cref{lem:c_continuous} does not apply because the Fréchet derivative with respect to $f$ loses derivatives. Instead, we extract the solution map $e(v) = f$ using the method of characteristics \cite[Chapter 2.2]{NicklLectureNotes2023}. 

Following \cref{rem:spillting_pde}, we lift the non-homogeneous boundary data. Let $\tilde{g} \in H^2(\mathcal{O})$ be the unique harmonic lift of $g$, satisfying $\Delta \tilde{g} = 0$ in $\mathcal{O}$ and $\tilde{g} = g$ on $\partial\mathcal{O}$. We expand the divergence in \cref{eq:darcy_standard} via the product rule and choose the linearised state variable $v = \Delta u$. Because $\Delta \tilde{g} = 0$, we have $v = \Delta(u - \tilde{g})$. 

Let $\mathcal{K} = \Delta^{-1} : L^2(\mathcal{O}) \to H^2(\mathcal{O}) \cap H^1_0(\mathcal{O})$ denote the inverse of the Dirichlet Laplacian. 
We can uniquely decompose the state as $u = \mathcal{K}v + \tilde{g}$. 
Substituting this into the expanded PDE transforms \cref{eq:darcy_standard} into a first-order linear transport equation for $f$, driven entirely by the velocity field induced by $v$ and the boundary lift $\tilde{g}$ i.e.
\begin{equation}\label{eq:darcy_ode}
    \nabla (\mathcal{K} v + \tilde{g}) \cdot \nabla f + v f = h.
\end{equation}

Because the characteristic curves depend on the unknown state $v$, the inflow boundary $\Gamma_{in}(v) := \{x \in \partial \mathcal{O} : \nabla (\mathcal{K} v(x) + \tilde{g}(x)) \cdot \nu(x) < 0\}$ fluctuates with $v$. 

\begin{assumption} \label{ass:darcy_structure}
To make the inverse map mathematically well-defined and stable, we impose the following structural conditions:
\begin{itemize}
    \item The ground truth variable $v_0$ satisfies a strict lower bound on its induced velocity field: $\inf_{x \in \mathcal{O}} |\nabla (\mathcal{K} v_0(x) + \tilde{g}(x))| \geq \delta > 0$. Furthermore, the field is non-trapping, meaning every characteristic curve exits the domain through the outflow boundary in finite time.
    \item The true parameter $f_0$ is perfectly known on a fixed, open boundary segment $\Gamma \subseteq \partial \mathcal{O}$. For $v_0$, the strict inflow boundary is compactly contained within this known segment: $\Gamma_{in}(v_0) \Subset \Gamma$. 
    \item The domain $\mathcal{O}$, boundary lift $\tilde{g}$, and source $h$ are sufficiently smooth such that the posterior concentrates on a sequence of sets $V_n \subset V$ where the vector fields $\nabla (\mathcal{K} v + \tilde{g})$ are uniformly bounded in $W^{1,\infty}(\mathcal{O})$ and the parameters $f = e(v)$ are uniformly bounded in $W^{1,\infty}(\mathcal{O})$.
\end{itemize}
\end{assumption}

\begin{lemma}
\label{lem:darcy_lip}
Suppose \cref{ass:darcy_structure} holds. For $v$ sufficiently close to $v_0$ in the $L^\infty$-norm, the fluctuating inflow boundary remains within the known data region: $\Gamma_{in}(v) \subset \Gamma$. Consequently, for any $v_1, v_2 \in V_n$, their respective parameters $f_1 = e(v_1)$ and $f_2 = e(v_2)$ share identical, fixed initial conditions $f_1 = f_2$ on $\Gamma_{in}(v_1)$. 
Furthermore, the solution map $e: V_n \to L^2(\mathcal{O})$ is Lipschitz continuous, satisfying:
\begin{equation}
    \|e(v_1) - e(v_2)\|_{L^2(\mathcal{O})} \leq \tilde{C} \|v_1 - v_2\|_{L^2(\mathcal{O})}, \quad \forall v_1, v_2 \in V_n.
\end{equation}
\end{lemma}
\begin{proof}
We bound the difference between $f_1 = e(v_1)$ and $f_2 = e(v_2)$ by subtracting their respective transport equations \cref{eq:darcy_ode}. Crucially, the fixed boundary lift $\tilde{g}$ cancels out in the difference of the velocity fields, $\nabla(\mathcal{K}v_1 + \tilde{g}) - \nabla(\mathcal{K}v_2 + \tilde{g}) = \nabla \mathcal{K}(v_1 - v_2)$, yielding:
\begin{equation}
    \nabla (\mathcal{K} v_1 + \tilde{g}) \cdot \nabla (f_1 - f_2) + (f_1 - f_2) v_1 = (v_2 - v_1) f_2 + \nabla \mathcal{K}(v_2 - v_1) \cdot \nabla f_2.
\end{equation}
We integrate this ordinary differential equation along the characteristic curves $X(t)$ defined by $\dot{X} = \nabla (\mathcal{K} v_1(X) + \tilde{g}(X))$. Because $\Gamma_{in}(v_1) \subset \Gamma$, the initial boundary error $(f_1 - f_2)$ is strictly zero at the inflow. Applying Grönwall's inequality along the flow lines provides the standard stability estimate (see \cite[Theorem 2.2.1]{NicklLectureNotes2023}):
\begin{equation}
    \|e(v_1) - e(v_2)\|_{L^2(\mathcal{O})} \leq C \left( \|v_1 - v_2\|_{L^2(\mathcal{O})} \|f_2\|_{L^\infty} + \|\nabla \mathcal{K}(v_1 - v_2)\|_{L^2(\mathcal{O})} \|\nabla f_2\|_{L^\infty} \right).
\end{equation}
Because $\mathcal{K} = \Delta^{-1}$ maps $L^2(\mathcal{O}) \to H^2(\mathcal{O}) \cap H^1_0(\mathcal{O})$, the gradient operator $\nabla \mathcal{K}$ is a bounded linear map from $L^2(\mathcal{O}) \to H^1(\mathcal{O})$. By the continuous embedding $H^1(\mathcal{O}) \hookrightarrow L^2(\mathcal{O})$, we have $\|\nabla \mathcal{K}(v_1 - v_2)\|_{L^2} \lesssim \|v_1 - v_2\|_{L^2}$. On the sets $V_n$, the norms $\|f_2\|_{W^{1,\infty}}$ are uniformly bounded by a constant independent of $n$, yielding the desired $L^2$-Lipschitz bound.
\end{proof}

\begin{corollary}
\label{cor:darcy_contraction}
Suppose that $v_0\in H_d^{\beta'}$, let $\tilde{v}_0 \in H_d^{\beta}$, and let \cref{ass:darcy_structure} hold. If the prior covariance $C_0$ and the stopping time $\tau_{\rm dp}$ are specified as in \cref{corollary:schrodinger_contraction_rate}, and if $\widetilde{\Pi}_n(V_n \mid \widetilde{Y}_n) \to 1$ in $P_{v_0}$-probability, then the posterior $\Pi_{\tau_{\rm dp},n}(f \mid Y_n)$ for the Darcy permeability field $f$ contracts at the optimal rate for $\tilde{v}_0$, $\epsilon_n \asymp n^{-2\beta/(2\beta + 2\alpha + 2p + 1 + d)}$, and its credible sets achieve asymptotic frequentist coverage $1$.
\end{corollary}
\begin{proof}
Because $e(v)$ is locally Lipschitz on $V_n$ by \cref{lem:darcy_lip}, it admits the linear modulus of continuity $\omega_n(t) = \tilde{C}t$. The rest of the proof follows identically as in \cref{corollary:schrodinger_contraction_rate} and \cref{corollary:schrodinger_coverage} by applying \cref{thrm:orginal_post_contracts} and \cref{corollary:coverage_original_post} directly apply, scaling the contraction rate and credible set radii by $\tilde{C}$.
\end{proof}

\section{Numerics and Ensemble Kalman Bucy Inversion} \label{sec:enkbf}
In this section, we provide a review of the Ensemble Kalman Bucy filter and the homotopy method, which allows us to see $\tau_n$ as a time parameter. We can then use this algorithm to verify the assumptions of \cref{sec:main_results}. To do this, we transfer the problem into sequence space. However, the theory is not limited to this; see \cite{tienstra2025, koers2024}—and thus the numerical implementation can be done in the discretised function space. We then describe how we can view $\tau_n$ in \cref{eq:loss_functional} as a time parameter that transforms an initial distribution into the target distribution. We also summarise the algorithmic details found in \cite{tienstra2025} and how to implement the early stopping of \cref{sec:main_results}. Furthermore,  in this section, we work in the $D(n)$ finite-dimensional spaces, and will denote the measures now by $\pi$. This theory holds in the infinite-dimensional setting. We first introduce the \textit{homotopy ansatz}
\begin{equation}
    \pi_\tau(\theta) \ \propto \  e^{-\frac{\tau}{2}(\mathcal{K}v - \widetilde{Y})^{\rm T}R^{-1}(\mathcal{K}v - \widetilde{Y})}\,
\pi_0(v)
\end{equation}
When $\tau=0$, we are in the prior. As $\tau\rightarrow1$, we reach the posterior by successive weighting of the likelihood. In \cite{tienstra2025}, they showed that this $\tau$ and the scale parameter of the prior covariance are related, in that setting $\tau = \tau_n$ and setting the prior covariance to the uncaled $C_0$ is the same as allowing $\tau=1$ and starting with prior covariance $\tau_n C_0$. Suppose then that we have samples, in the setting of a filter called particles from $\pi_0$. What we want is then some algorithm that iteratively updates the particles such that at the discrete time $t$, the particles are approximately samples from $\pi_t$.  If $\mathcal{K}$ is linear, all $\pi_\tau$ are normal with mean and covariance given by 
\begin{align}
m_\tau &= m_0 - C_0 \mathcal{K}^{\rm T} (\mathcal{K}  C_0^{(\gamma)} G^{\rm T} + \tau^{-1} R)^{-1} (Gm_0 - y),\\
C_\tau &= C_0 - C_0 \mathcal{K}^{\rm T} (\mathcal{K} C_0 
\mathcal{K}^{\rm T} + \tau^{-1} R)^{-1} \mathcal{K} C_0
\end{align}
Consider then then, Kalman--Bucy mean field filter equations (EnKBF) \cite{reich10}: 
\begin{equation} \label{eq:EnKBF}
{\rm d}\mathcal{V}\tau = C_\tau \mathcal{K}^{\rm T} R^{-1} \left\{ (y-G\mathcal{V}_\tau {\rm d}\tau -
R^{1/2} {\rm d}W_\tau \right\}
\end{equation}
with initial conditions drawn from the prior, that is, $\mathcal{V}_0^{(\gamma)} \sim \mathcal{N}(m_0,\gamma C_0)$. 
Here $W_\tau$ denotes standard $d_y$-dimensional Brownian motion. We then have that 
\begin{equation}
\mathcal{V} \sim \pi_\tau
\end{equation}
for all $\tau>0$. 
The discrete time formulations of \cref{eq:EnKBF} can be written as 
 
\begin{align} 
\label{eq:DEnKBFBF}
\frac{\rm d}{{\rm d}t}\mathcal{V}_t &= -\frac{1}{2} n \Sigma^J_t \mathcal{K}^{\rm T}  \left(\mathcal{K}\mathcal{V}_t+\mathcal{K}m_t-2y \right) .
\end{align}
We can then implement the ENKF as follows. Let the number of particles of the ensemble be denoted by $J$. The discrete time index will now be denoted by $t_k \ge 0$. The ith particle of the ensemble of $J$ particles at time $t_k$  by $v_k^{(i)}$. The initial condition is chosen as 
\begin{equation} \label{eq:ICs}
v_0^{(i)} \sim {\rm N}(0,C_0)
\end{equation}
for $i=1,\ldots,J$.  The empirical mean of the ensemble at time $t_k$ is denoted by 
\begin{equation}
   m_k^{J}= \frac{1}{J}
    \sum_{i=1}^J v_k^{(i)}
\end{equation}

Similarly, the empirical covariance matrices are given by $C_k^J$. To compute the Kalman gain, we need to introduce the empirical covariance matrix
between $\mathcal{V}$ and $\mathcal{K}\mathcal{V}$, which we will denote by $\mathcal{C}_n^{J} \in \mathbb{R}^{D(n) \times D(n)}$, 
as well as the empirical covariance matrix of  $\mathcal{K}\mathcal{V}$, which will denoted by $S_k^{J} \in \mathbb{R}^{D(n)\times D(n)}$. For clarity, they are given below:
\begin{equation}
\label{eq:em_cov}
    S_k^J = \frac{1}{J-1} \sum_{i=1}^J (
    \mathcal{K}v_k^{(i)} -  m^J_{\mathcal{K},k})(m^J_{\mathcal{K},k} - \mathcal{K}v_k^{(i)})^{\rm T},
\end{equation}
where $m^J_{\mathcal{K},k}$ denotes the empirical mean of $\mathcal{K}v$. Similarly,
\begin{equation}
\label{eq:em_cross_cov}
    \mathcal{C}_k^{J} = \frac{1}{J-1} \sum_{i=1}^J (
    \theta_k^{(i)} - m^J_k)(Gv_k^{(i)}) - m^J_{\mathcal{K},k})^{\rm T}.
\end{equation}

The deterministic discrete time update formulas, which approximate \cref{eq:EnKBF}, are then
\begin{subequations} \label{eq:numerical_DEnKBF}
\begin{align}
v_{{k+1}}^{(i)} &= v_{k}^{(i)} - \frac{1}{2} K_{k} \left(\mathcal{K}v_{k}^{(i)} +m^J_{\mathcal{K},k} - 2y
\right)
\end{align}
\end{subequations}
where the ith Kalman gain matrix is 
\begin{equation} \label{eq:Kalman_gain_numerical}
K_{k} = \Delta t \,\mathcal{C}_{k}^J \left(\Delta t S_{k}^J + I \right)^{-1}.
\end{equation}

The standard discrepancy principle stops the iteration of the EnKBF  whenever 
\begin{equation} \label{eq:numerical_dp_stopping}
   k_{\rm dp} = \inf \left\{k \geq k_0 : \|G(\tilde{v}^J_k) - Y\|^2 \leq \kappa \right\}.
\end{equation}
and we then chose $\kappa = C D(n) / n$ and $0<C \le 1$ and $k_0$ is the initial time. We then have the resulting algorithm which implements \cref{eq:DEnKBFBF} and \cref{eq:dp_stopping_time}.

\newenvironment{algocolor}{%
   \setlength{\parindent}{0pt}
   \itshape
}{}

\begin{algorithm} 
      \caption{Deterministic EnKF}
      \label{alg:deter_enkf}
      \begin{algocolor}
      \begin{algorithmic}
      \Require $J > 0,\  m_0, \ C_0, \ y,\ \mathcal{K}$ 
       \State $V_0 \gets \texttt{initialize}(J, m_0, C_0)$ \Comment{$\Theta  \in \mathbb{R}^{D(n) \times J}$}
        \State $R_0 \gets ||\mathcal{K}( m^J_0) - y ||^2 $ 
        \State $\kappa_{\rm dp} = D(n)*n$ \Comment{see \eqref{eq:dp_stopping_time}}
       \While{$R_k < \kappa_{\rm dp}$}  
            \State $K_{k} \gets \Delta t \,\mathcal{C}_{k}^J \left(\Delta t S_{k}^J + I \right)^{-1}$ \Comment{$\mathcal{C}_{k}^J$ see \eqref{eq:em_cross_cov},  $\Sigma_{k}^J$ see \eqref{eq:em_cov}}
           \For {$i \in \{1,..,J\}$} 
           \begin{align*}
               v_{{k+1}}^{(i)} &= v_{k}^{(i)} - \frac{1}{2} K_{k} \left(\mathcal{K}v_{k}^{(i)} +m^J_{\mathcal{K},k} - 2y
            \right)
           \end{align*}
           \EndFor
          \State $R_{k+1} \gets ||\mathcal{K }\left(m^J_{k+1}\right) - y||^2 $ 
       \EndWhile
       \State \textbf{Return} $V_k$
    \end{algorithmic}   
  \end{algocolor}
\end{algorithm}

%%%%%%%%%%%%%%%%%%%%%%%%%%%%%%%%%%%%%%

\subsection{Numerical Results for Schrödinger} \label{subsec:numerical_results} 
\begin{figure}[h!]
    \centering
    \includegraphics[width=1\textwidth]{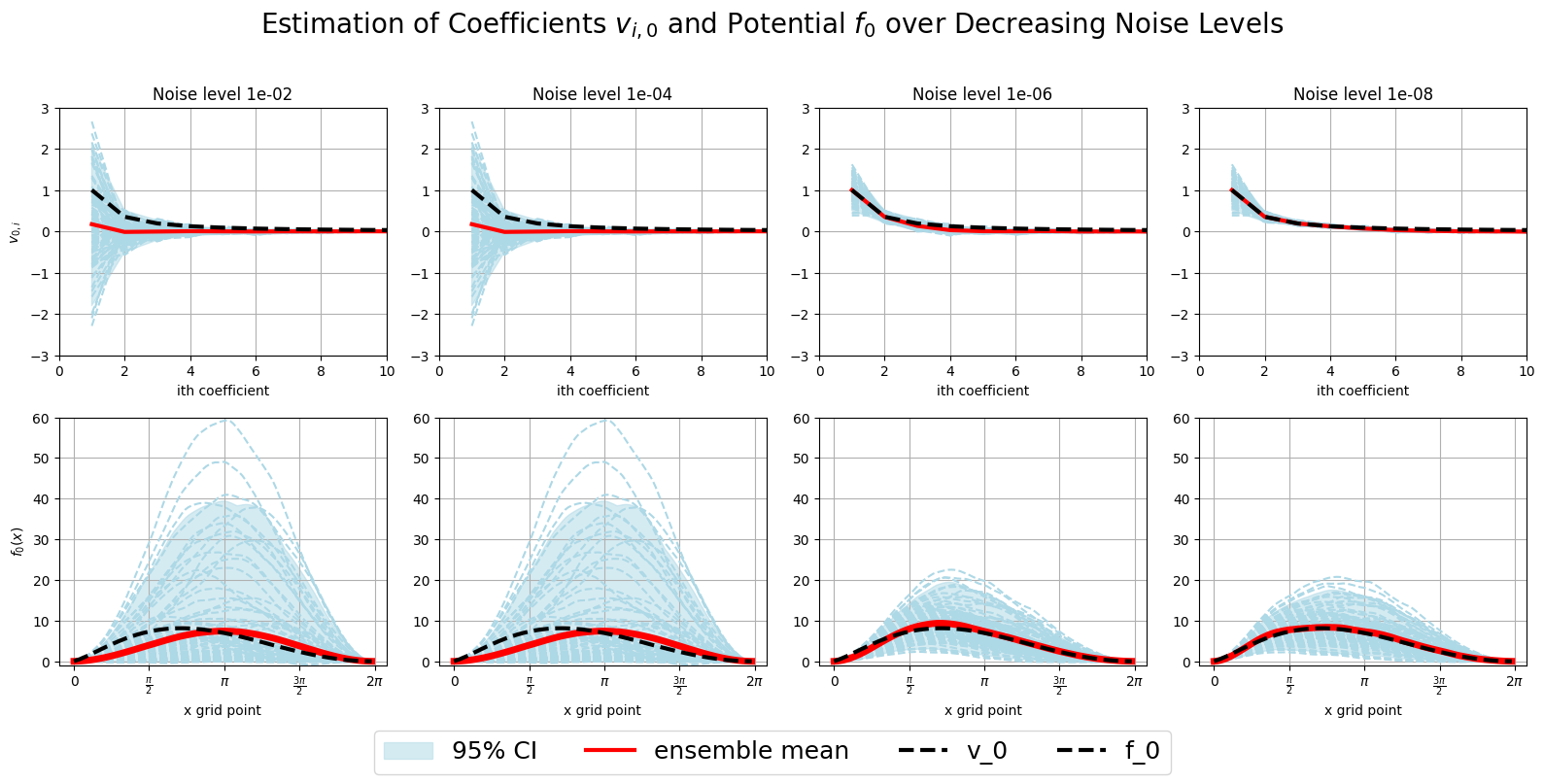}
    \caption{Here we plot the results of running EKI with early stopping, \cref{alg:deter_enkf}, on the Schrödinger problem. From left to right, the noise level decreases. On the top is the estimation for $v_{0,i}$, the coefficients of $v_0$. We plot only the first $10$ coeffients as the remaining are essentially zero, and this zoomed-in perspective shows how uncertainty in the coefficients propagates to the uncertainty in $f_0$.  On the bottom are the resulting transformed estimates for $f_0$ in function space over a grid of $100$ points,  and  $\kappa$ is chosen to be $D(n)*n$ where $D(n)=100$. For all noise levels, we fix the grid and only scale the variance of the noise in the linear observations. The red solid line is the ensemble mean, the black dashed line is the ground truth, the blue thinner dashed lines are the ensemble particles, and the blue filled region is the $95\%$ credible region computed by taking the  $95\%$ quantiles of the ensemble. }
    \label{fig:schrodinger}
\end{figure}

We demonstrate the results \footnote{\url{https://github.com/Tienstra/EarlyStoppingBIP}} of \cref{sec:main_results} for the Schrödinger equation in 1-dimension on $[0,1]$. We chose $v_0= \mathcal{L}u_{f_0}$, the ground truth, to be such that $v_0(x) = \sum_i v_{i,0} \phi_i(x)$ where $v_{i,0} = i^{-5/2}$ and $\phi_i(x) = \sqrt{2}\sin(i\pi x)$ are the Dirichlet eigen functions of $\mathcal{K}$, the inverse  of the negative Laplace operator.  We chose homogenous boundary conditions $g(0)=0$ and $g(1) = 0$ and so $\tilde{g} = g$.  The covariance operator has eigenvales $\lambda_i = i^{-1/2-\alpha}$ where $\alpha =2$. We then ran EKI, see  \cref{alg:deter_enkf} with early stopping rule \cref{eq:numerical_dp_stopping} with $C=1$ in sequence space to recover a finite number $D(n)$ of coefficients of $v_0$. We then transformed the estimates back into the function space, via this basis, to have an estimator for $v_0$, which we denote by $v_{\tau_{\rm dp}}$. As $v = \mathbb{L}u_{f}$ and every $u_f$ is uniquely determined by $f$, we thus have an estimator for $f_0$ via the solution map \cref{eq:solution_map_e}. That is
\begin{equation}
    f = e(\mathbb{L}u_f) = \frac{v}{2(\mathcal{K}v + \tilde{g})}
\end{equation}
when $\rm{essinf } \ \mathcal{K}v + \tilde{g} >0$ and zero else. We remark that when $f$ is positive, a unique solution to \cref{eq:generic_pde} is guaranteed. As $\mathcal{K} = \mathbb{L}^{-1} = - \Delta^{-1}$, and the eigenvalues of $\mathcal{L}$ are positive \cref{eq:decay_laplace_eigenvalues} so are the eigenvalues of $\mathcal{K}$ and so $v$, and $\mathcal{K}v$ will have the same sign. That ensures that $f$ will always be positive. The results can be found in \cref{fig:schrodinger}.
We can see in \cref{fig:schrodinger}, that as $v_{\tau_{\rm dp}} \rightarrow v_0$, that $f_{\tau_{\rm dp}}$ also converges to $f_0$. We also see that the posterior spread estimated by the ensemble spread shrinks as $n \rightarrow \infty$ and that $f_0$ is within the $95\%$ quantiles of the particles.  While the results of \cref{sec:main_results} are asymptotic, one can ask if this method works in the finite $n$ setting. We see from the \cref{fig:schrodinger} that the answer seems positive as the finite dimension was fixed at $100$. 

\section{Conclusion} \label{sec:conclusion}

In this paper, we developed a methodology for tuning Gaussian priors for non-linear Bayesian inverse problems arising from semilinear PDEs. Our approach builds on the transformation technique introduced in \cite{koers2024}, which reformulates the original non-linear estimation problem for $f_0$ as a linear problem for $v_0 = \mathcal{L}u_{f_0}$. We extended this linearisation method to a wide class of semilinear elliptic PDEs, and proposed an extension to time-dependent parabolic equations.

We also extended the linearisation method to more general solution maps $e$, see \cref{eq:solution_map_e}, which admit a modulus of continuity. With this more general notion of continuity, we proved in \cref{proposition:transfer} and \cref{proposition:credible} that rates of posterior contraction and credible-set coverage can be transferred from the linearised posterior to the pulled-back posterior for $f$. We further provided a condition on $c(u_f,f)$, see \cref{prop:nemytskii_modulus}, under which such solution maps $e$ exist.

Using the linearised setting, the posterior distribution for $v_0$ is Gaussian, allowing us to apply early stopping, guided by the discrepancy principle of \cite{tienstra2025}, to determine the scale parameter of the prior covariance. As shown in \cite{tienstra2025}, this procedure yields a posterior that contracts at a near-optimal rate for the linearised problem. Via the solution map $e$, these contraction properties are then transferred to the posterior ${\Pi}_n(f \mid \widetilde{Y})$ for the original parameter $f_0$, with the rate modified by the modulus of continuity when the map is not Lipschitz. Analogous results hold for credible-set coverage, up to a corresponding change in the radius.

To demonstrate our general results, we first analysed the canonical example of the stationary Schr\"odinger equation, for which we prove posterior contraction and frequentist coverage and also provide numerical experiments. We then showed how the same framework applies to the stationary Allen-Cahn equation, where the polynomial non-linearity fits naturally into our general condition on $c(u_f,f)$. We also treated the parabolic Allen-Cahn equation as a concrete example of the time-dependent extension. Finally, we discussed Darcy flow as an edge case: although it is not semilinear, the problem can still be brought into a form where a solution map may be constructed under additional structural assumptions.

We also showed how the prior scale parameter can be interpreted as the homotopy parameter in \cite{reich10}, and identified it with time in the ensemble Kalman--Bucy filter. This gives an implementable algorithm in which early stopping becomes a stopping time for the filter. The numerical results for the stationary Schr\"odinger equation support the asymptotic theory in the finite-dimensional setting.

We note that the theoretical development in this work, as well as in \cite{koers2024}, relies on the existence of an operator $\mathcal{K}$ that enables the reformulation of the original problem into the structure given in \cref{eq:generic_pde}. An interesting direction for future research would be to investigate whether similar techniques can be extended to settings where $\mathcal{K}$ serves only as an approximation that linearises the problem, rather than as an exact inverse. This would broaden the applicability of the method to non-linear problems not directly associated with linear differential operators.

\appendix
% make \cref say “Appendix” instead of “Section” *inside* the appendix
\crefalias{section}{appendix}

%%%%%%%%%%%%%%%%%%%%%%%%%%%%%%%%%%%%%%

\section{Appendix A} \label{sec:appendix_A}

%%%%%%%%%%%%%%%%%%%%%%%%%%%%%%%%%%%%%%
In this section we collect important definitions refered to and  were place here to maintain redabilty and focus of the manuscript.

\begin{definition}[Emperical Inner Product and Emperical Projection Operator]
    Write $\mathcal{H}:=L^{2}([0,1]^{d})$ with inner product $\langle\cdot,\cdot\rangle$ and
    orthonormal basis $(\psi_{i})_{i\in\mathbb{N}}$, and set
    $S_{n}:=\operatorname{span}(\psi_{1},\dots,\psi_{n})\subset\mathcal{H}$.
    Given the design $X_{1},\dots,X_{n}$ define the empirical inner product
    \[
    \langle g,h\rangle_{n}:=\frac{1}{n}\sum_{j=1}^{n}g(X_{j})h(X_{j}),
    \qquad
    \langle \mathbf{Y},h\rangle_{n}:=\frac{1}{n}\sum_{j=1}^{n}Y_{j}\,h(X_{j}),
    \]
    the design matrix $\Psi_{n}:=\big(\psi_{i}(X_{j})\big)_{1\le i,j\le n}\in\mathbb{R}^{n\times n}$,
    and the empirical projection operator
    \begin{equation}
    \label{eq:empirical_projection_operator}
    \Pi_{n}:\mathcal{H}\to S_{n},
    \qquad
    \Pi_{n}g:=\sum_{i=1}^{n}\langle g,\psi_{i}\rangle_{n}\,\psi_{i}.
    \end{equation}
    In the basis $(\psi_{i})_{i\le n}$ the restriction $\Pi_{n}\!\mid_{S_{n}}:S_{n}\to S_{n}$ has matrix
    $G_{n}=\big(\langle\psi_{i},\psi_{k}\rangle_{n}\big)_{i,k\le n}=n^{-1}\Psi_{n}\Psi_{n}^{\top}$,
    which is symmetric positive semi-definite with $\mathbb{E}[G_{n}]=\mathrm{Id}_{S_{n}}$;
    hence $(\Pi_{n}\!\mid_{S_{n}})^{1/2}$ is well defined, and it is invertible iff $\Psi_{n}$ is.
\end{definition}

\section{Appendix B} \label{sec:appendix_B}

In this section, we repeat the necessary results from \cite{koers2024} and \cite{tienstra2025}, respectively. The original statements and proofs can be found in the sources.

%%%%%%%%%%%%%%%%%%%%%%%%%%%%%%%%%%%%%%
\subsection{Results for Linear to Non-linear}
%%%%%%%%%%%%%%%%%%%%%%%%%%%%%%%%%%%%%%

\begin{proposition} (Proposition 2.1 in \cite{koers2024})
\label{proposition:koers_contraction_pi_tilde}
Suppose that the posterior distribution 
$$ \tilde{\Pi}_n(v \in \cdot \mid \widetilde{Y}_n) $$ 
of $v$ in model (L) contracts under $v_0$ to $v_0 = \mathcal{L} u_{f_0}$ at rate $\epsilon_n$ in $(V, \|\cdot\|)$ and satisfies 
$$ \Pi_n(v \in V_n \mid \widetilde{Y}_n) \xrightarrow{P} 1 $$
for given sets $V_n \subset V$. If \cref{eq:solution_map_e} holds for a map $e$ such that $e: V_n \to L_2$ is Lipschitz at $v_0$, then the posterior distribution of $f$ in model (N), \cref{eq:nonlinear_model}, attains a rate of contraction $\epsilon_n$ under $f_0$ relative to the $L_2$-norm.
\end{proposition}

\begin{proposition}(Proposition 2.2 in \cite{koers2024})
    \label{proposition:koers_credible_sets_equal}
    The credible levels of credible ses $\widetilde{C}_n(\widetilde{Y}_n)$ for $v$ in \cref{eq:linearized_post} and $C_n(Y_n)$ for $f$ in \cref{eq:original_post} are equal, and so are the coverage levels of the sets at $v_0=\mathcal{L}u_{f_0}$ and $f_0$ respectively. Furthermore if the map $e:V_n \rightarrow L_2$ in \cref{eq:solution_map_e} is uniformly Lipschitz at points $\bar{v}_n \in V_n$ and  $\widetilde{C}_n(\widetilde{Y}_n) \subseteq V_n$, then on the event $\bar{v}_n \in  \widetilde{C}_n(\widetilde{Y}_n)$ the $L_2$ diameter of the sets $C_n(Y_n)$ under $f_0$ are of the same order in probability as the $\|\cdot \|-$diameters of the set $\widetilde{C}_n(\widetilde{Y}_n)$ under $v_0$.
\end{proposition}

%%%%%%%%%%%%%%%%%%%%%%%%%%%%%%%%%%%%%%
\subsection{Results for data-driven posterior}
%%%%%%%%%%%%%%%%%%%%%%%%%%%%%%%%%%%%%%

\begin{theorem} (Theorem 2.1 in \cite{tienstra2025})
\label{thrm:tien_contraction_rate_early_stopping}
    Let $v_0 \in H_d^{\beta^\prime}$, and denote the posterior associated to the estimated stopping time  $\tau_{\rm dp }$ by $\widetilde{\Pi}_{n,\tau_{\rm dp }}(\cdot \mid \widetilde{Y}_n)$. Let $\mathcal{K}$ and $C_0$ have the same eigenfunctions. Define $\tilde{v}_0 :=C_0^{-1/2}v_0 \in H_d^\beta$. Denote the eigenvalues of $\mathcal{K}$ and $C_0$ by $\sigma_i$ and $\lambda_i$ respectively. If the eigenvalues of $\mathcal{K}$ have polynomial decay, that is
    \begin{equation}
        \sigma_i \asymp i^{-p/d}
    \end{equation}
    And we chose entry-wise prior
    \begin{equation}
        v_i \sim \mathcal{N}(0,\tau_n^2 \lambda_i). \quad \lambda_i \asymp i^{-1-2\alpha/d}.
    \end{equation}
    such that $\beta < 2 \alpha + 2p + d$,
     then
    \begin{equation}
       \widetilde{\Pi}_{n,\tau_{\rm dp }} \left( \widehat{v}_{\tau_{\rm dp }} : || \widehat{v}_{\tau_{\rm dp }} - v_0||^2_{\ell^2(\mathbb{N})} \geq M_n \epsilon_n\right) \rightarrow 0
    \end{equation}
    for every $M_n \rightarrow \infty$, and with $\epsilon_n \asymp n^{-2\beta / (2\beta + 2p + 2\alpha + 1 + d)}$.
\end{theorem}

\begin{corollary}
\label{corollary:tien_assympotic_coverage}
    For fixed $\alpha > 0$, if $v_{0,i} = C_i i^{-1-2\beta'/d}$ for all $i=1,...,D(n)$, and $\beta \leq d + 2\alpha +2p$ where $\widetilde{v}_0 \in H_d^\beta$ then as $n \rightarrow \infty$, $ \widetilde{\Pi}_{n, \tau_{\rm dp}}$  has frequentist coverage 1.
\end{corollary}

\section{Appendix C} \label{sec:appendix_C}
In this section we present slight modifications of established results in PDE theory and functional analysis.

\begin{proof}[Proof of \cref{prop:nemytskii_modulus}]
The following proof utilizes standard bounds for Nemytskii operators in Lebesgue spaces, which rely heavily on the interplay between Carathéodory growth conditions and Hölder's inequality (see, for instance \cite{zeidler_nonlinear_1990}). 

Let $B \subset U$ and $B_f \subset \mathcal{F}$ be bounded sets. 
Because $U \subset H^{2m}(\mathcal{O})$ and $\mathcal{F} \subset H^l(\mathcal{O})$, the Sobolev embedding theorem implies that for every $\alpha \in A_{\mathrm{loc}}(k, 2m)$ and $\beta \in A_{\mathrm{loc}}(l, l)$, 
\[
\sup_{u \in B}\|D^\alpha u\|_{L^\infty(\mathcal{O})} < \infty, \qquad \sup_{f \in B_f}\|D^\beta f\|_{L^\infty(\mathcal{O})} < \infty.
\]
Furthermore, for every $\alpha \in A_{\mathrm{glob}}(k, 2m)$ and $\beta \in A_{\mathrm{glob}}(l, l)$, the choices of exponents $q_\alpha$ and $p_\beta$ in (C1) guarantee the continuous embeddings $U \hookrightarrow L^{2q_\alpha}(\mathcal{O})$ and $\mathcal{F} \hookrightarrow L^{2p_\beta}(\mathcal{O})$. 
Therefore, we can choose constants $R, R_f > 0$ uniformly bounding the local jets almost everywhere.
By (C1), taking the $L^2$-norm and applying the continuous embeddings yields
\[
\|C(u;f)\|_V \leq \|a_{R,R_f}\|_{L^2} + C_u \sum_{\alpha \in A_{\mathrm{glob}}(k, 2m)} \|u\|_U^{q_\alpha} + C_f \sum_{\alpha \in A_{\mathrm{glob}}(l, l)} \|f\|_{\mathcal{F}}^{p_\alpha}.
\]
Since $\|u\|_U$ and $\|f\|_{\mathcal{F}}$ are uniformly bounded on $B$ and $B_f$, we conclude $\sup_{u\in B, f\in B_f} \|C(u; f)\|_{V} < \infty$ and confirm $C$ maps bounded sets into bounded sets.
To establish the modulus of continuity with respect to $u$, fix $f_0 \in \mathcal{F}$ and choose $r,s > 0$ such that $B_r(u_{f_0})$ and $B_s(f_0)$ are bounded. 
Next, we apply (C2) to the difference $C(u; f)(x) - C(v; f)(x)$. 
As the parameter $f$ is fixed in this evaluation, the parameter jets are identical i.e. $\zeta_1 = \zeta_2$. 
Therefore, all terms in (C2) involving differences of $\zeta$ evaluate identically to zero and we obtain
\begin{align*}
|C(u; f)(x) - C(v; f)(x)| 
\leq & L_{R,R_f} \omega_{R,R_f}(|J^k_{\mathrm{loc}}(u-v)(x)|) 
\\
&+ L_{R,R_f} \sum_{\alpha\in A_{\mathrm{glob}}(k, 2m)} ( |D^\alpha u(x)|^{q_\alpha-\gamma_\alpha} + |D^\alpha v(x)|^{q_\alpha-\gamma_\alpha} ) |D^\alpha(u-v)(x)|^{\gamma_\alpha} .
\end{align*}

For the local term, because $t \mapsto \omega_{R,R_f}(\sqrt{t})^2$ is concave, Jensen's inequality gives
\[
\| \omega_{R,R_f} (|J^k_{\mathrm{loc}}(u-v)|) \|_{L^2}
\leq 
|\mathcal{O}|^{1/2} \omega_{R,R_f} ( |\mathcal{O}|^{-1/2} \|J^k_{\mathrm{loc}}(u-v)\|_{L^2} ).
\]
Since the jet mapping is bounded from $U$ to $L^2$, there exist constants $w_1, w_2 > 0$ such that this is further bounded by $w_2 \omega_{R,R_f}(w_1\|u-v\|_U)$.

For the global terms i.e. $\alpha \in A_{\mathrm{glob}}$, we apply Hölder's inequality with conjugate exponents $\frac{2q_\alpha}{q_\alpha - \gamma_\alpha}$ and $\frac{2q_\alpha}{\gamma_\alpha}$ and obtain
\begin{align*}
\| ( |D^\alpha u|^{q_\alpha-\gamma_\alpha} + |D^\alpha v|^{q_\alpha-\gamma_\alpha} ) |D^\alpha(u-v)|^{\gamma_\alpha} \|_{L^2}.
&\leq 
( \|D^\alpha u\|_{L^{2q_\alpha}}^{q_\alpha-\gamma_\alpha} + \|D^\alpha v\|_{L^{2q_\alpha}}^{q_\alpha-\gamma_\alpha} ) \|D^\alpha(u-v)\|_{L^{2q_\alpha}}^{\gamma_\alpha}.
\end{align*}
Since $u,v \in B_r(u_{f_0})$, the first factor is bounded by a constant dependent on $r$. 
The second factor is bounded by $c_\alpha \|u-v\|_U^{\gamma_\alpha}$ due to the continuous embedding.

Combining the local and global estimates, there exists a constant $w_0 > 0$ such that
\[
\|C(u; f) - C(v; f)\|_V
\leq 
w_0 ( \omega_{R,R_f}(w_1\|u-v\|_U) + \sum_{\alpha\in A_{\mathrm{glob}}} \|u-v\|_U^{\gamma_\alpha} ).
\]
This establishes \cref{eq:C_modulus} with the explicit modulus $\omega_C$ from \cref{eq:explicit_C_modulus}. 
Finally, because $\omega_{R,R_f}(t) \to 0$ as $t \downarrow 0$ and $\gamma_\alpha > 0$, we conclude that $\omega_C(t) \to 0$ as $t \downarrow 0$.
\end{proof}

\begin{proof}[Proof of \cref{rem:frechet_derivative}]
  We choose $A = \partial_f C(u_{f_0}; f_0)$.
  Let $f_1, f_2 \in B_s(f_0)$ and define the convex combination $f_t = t f_1 + (1-t) f_2$. 
  As open balls in normed spaces are convex, we conclude $f_t \in B_s(f_0)$ for every $t \in [0, 1]$. 
  Next, we define the curve $\gamma \colon [0, 1] \to V$ by
  \[
    \gamma(t) := C(u; f_t) - t A (f_1 - f_2)
  \]
  for a fixed $u \in B_r(u_{f_0})$. 
  As the operator is partially Fréchet differentiable on $B_r(u_{f_0}) \times B_s(f_0)$, we apply the chain rule to obtain
  \[
    \gamma'(t) = \partial_f C(u; f_t) (f_1 - f_2) - A(f_1 - f_2)
  \]
  for $t \in (0,1)$. 
  Subsequently the mean value theorem for vector-valued functions gives
  \begin{align*}
    \|C(u; f_1) - C(u; f_2) - A(f_1 - f_2)\|_V &= \|\gamma(1) - \gamma(0)\|_V \\
    &\leq \sup_{t\in (0,1)} \|\gamma'(t)\|_V \\
    &\leq \sup_{t\in (0,1)} \| \partial_f C(u; f_t) - \partial_f C(u_{f_0}; f_0)\|_{L(\mathcal{F}, V)} \|f_1 - f_2\|_{\mathcal{F}}.
  \end{align*}
  Since the map $(u, f) \mapsto \partial_f C(u; f)$ is continuous at $(u_{f_0}, f_0)$, we can choose the radii $r$ and $s$ small enough such that 
  \[
    \sup_{t \in (0,1)} \| \partial_f C(u; f_t) - \partial_f C(u_{f_0}; f_0)\|_{L(\mathcal{F}, V)} \leq \frac{\theta}{\|A^{-1}\|}
  \]
  for all $u \in B_r(u_{f_0})$ and $f_1, f_2 \in B_s(f_0)$. This establishes the bound in \cref{eq:param_approx}.
\end{proof}

\begin{theorem}[Implicit Function Theorem \cite{hildebrandt1927}]
  \label{thm:graves}
  Let $X$, $Y$ be Banach spaces. 
  Let $x_0, \in X, y_0\in Y$ and let $B_r(x_0),B_s(y_0)$ be open balls around $x_0,y_0$ and radius $r,s$ respectively. 
  Let $G:B_r(x_0)\times B_s(y_0)\rightarrow Y$ with $G(x_0,y_0)=0$ and such that there exist a bounded linear operator $A\in L(Y,Y)$ with bounded inverse, 
  $0<\gamma<1/\|A^{-1}\|$ and a continuity modulus $w$ such that the conditions 
  \begin{align}
    \label{eq:HG_first_variable}
    \| G(x_1,y) - G(x_2,y)\|_Y &\leq w(\|x_1-x_2\|_X),\\
    \label{eq:HG_second_variable}
    \| G(x,y_1) - G(x,y_2) - A (y_1-y_2)\|_Y &\leq \gamma \|y_1-y_2\|_Y
%    \|G(x,y_0)\|_Y &\leq \frac{1-\gamma \|A^{-1}\|}{\|A^{-1}\|} \delta
  \end{align}
  hold true for every $x \in B_r(x_0)$, $y_1,y_2 \in B_s(y_0)$.

  Then, there exists $\tilde{r} \in (0,r)$ and a map $\phi: B_{\tilde{r}}(x_0) \rightarrow B_{s}(y_0)$ satisfying $G(x,\phi(x)) = 0$ for every $x\in B_{\tilde{r}}(x_0)$.
\end{theorem}
\begin{proof}
  We apply Theorem 3 in the original reference \cite{hildebrandt1927}.
  To this end we check the preconditions.
  First, we remark that Banach spaces are complete metric spaces and for linear operators the modulus is the operator norm.
  Next, we see that our assumptions should imply $y_*=0$ as we satisfy the equation $G(x_0,y_0)$ at least on the reference ball.
  However we still must check condition (H3) \cite[Theorem 3]{hildebrandt1927}.
  %To this end we have to prove existence of constant $\delta \in (0,s)$ such that
  %\[
  %  \|G(x,y_0)\|_Y \leq \frac{1-\gamma \|A^{-1}\|}{\|A^{-1}\|} \delta
  %\]
  %holds true for every $y_0 \in B_s^Y$ and all suiteable $x$.
  A continuity modulus is non-decreasing $\lim_{t \downarrow 0} w(t)=0$. 
  Therefore we find 
  $\tilde{r} \in (0,r)$ for every choice of $\delta \in (0,s)$ such that
  \[ 
    w(\tilde{r}) \leq \frac{1-\gamma \|A^{-1}\|}{\|A^{-1}\|} \delta 
  \]

  $G(x,y_0)$ satisfyies a continuity modulus $w$ on $B_r(x_0)$ and therefore on $B_{\tilde{r}}(x_0) \subseteq B_r(x_0)$.
  Therefore we conclude $G(x,y_0)\leq w(\tilde{r})\leq \frac{1-\gamma \|A^{-1}\|}{\|A^{-1}\|} \delta$ for all $x\in B_{\tilde{r}}(x_0)$.
  The statment now directly follows from \cite{hildebrandt1927}
\end{proof}
\begin{corollary}
  \label{cor:graves}
  Under the assumptions of \cref{thm:graves} we additionally 
  the implicitly given map $\phi$ satisfies the following scaled modulus of continuity  
  \begin{equation}
    \|\phi(x_1)-\phi(x_2)\|\leq \frac{\|A^{-1}\|}{1-\gamma \|A^{-1}\|} w(\|x_1-x_2\|)
  \end{equation}
  for every $x_1,x_2 \in B_{\tilde{r}}(x_0)$.
\end{corollary}
\begin{proof}
  For every $x_1,x_2\in X_0$ we have $G(x_1,\phi(x_1)) =0$. 
  Therefore, and due to \cref{eq:HG_first_variable,eq:HG_second_variable} in the assumption of the theorem we have 
  \[
    \begin{aligned}
      &\| \phi(x_1) - \phi(x_2)\|_Y \\
      &= \| \phi(x_1) - \phi(x_2) - A^{-1} (G(x_1,\phi(x_1)) - G(x_2,\phi(x_2)))\|_Y \\
      &\leq \|A^{-1}\|\| A (\phi(x_1)-\phi(x_2)) - G(x_1,\phi(x_1) + G(x_1,\phi(x_2)) \|_Y \\
      &\phantom{=}+ \|A^{-1}\| \| G(x_2,\phi(x_2)))- G(x_1,\phi(x_2)) \|_Y\\
      &\leq \|A^{-1}\| (\gamma \|\phi(x_1)-\phi(x_2)\|_Y + w(\|x_1-x_2\|_X))
    \end{aligned}
  \]
  Subtracting $\gamma \|A^{-1}\| \|\phi(x_1)-\phi(x_2)\|_Y$, we obtain the statement as $1-\gamma \|A^{-1}\|$ is positive.
\end{proof}

\paragraph{Acknowledgements.}
This work has been funded by Deutsche Forschungsgemeinschaft (DFG) - Project-ID 318763901 - SFB1294. 
Furthermore, G.H. has been supported by the European
Union under the Horizon Europe research and innovation programme (grant agreement no. 101188131, UrbanAIR).
The authors thank Sebastian Reich for his insightful feedback on multiple drafts of this manuscript. 
The authors also thank Botond Szabó for pointing out the connection between \cite{koers2024} and \cite{tienstra2025}. Finally the authors the reviewers of this paper for their constructive and insightful input. 
\newpage
%%%%%%%%%%%%%%%%%%%%%%%%%%%%%%%%%%%%%%
\nocite{*}
\bibliographystyle{abbrvnat}
%%%%%%%%%%%%%%%%%%%%%%%%%%%%%%%%%
%
\bibliography{refs}
%
%%%%%%%%%%%%%%%%%%%%%%%%%%%%%%%%
\end{document}